\documentclass[final,3p,times]{elsarticle}

\usepackage{amssymb}
\usepackage{graphicx,setspace}
\usepackage{psfrag}
\usepackage{amsfonts}
\expandafter\let\csname equation*\endcsname\relax
\expandafter\let\csname endequation*\endcsname\relax
\usepackage{amsmath}
\usepackage{amssymb, amsthm}
\usepackage{mathrsfs}
\usepackage{epstopdf}
\usepackage{float}
\usepackage{lineno,hyperref}
\usepackage{color}
\usepackage{subfigure}

\usepackage{amsmath}
\usepackage{dsfont} 
\usepackage{amssymb} 
\usepackage{graphicx,color}
\usepackage{extarrows}

\usepackage{graphicx}

\usepackage{epstopdf}

\journal{arXiv}

\begin{document}
\newtheorem{definition}{Definition}[section]
\newtheorem{lemma}{Lemma}[section]
\newtheorem{remark}{Remark}[section]
\newtheorem{theorem}{Theorem}[section]
\newtheorem{proposition}{Proposition}
\newtheorem{assumption}{Assumption}
\newtheorem{example}{Example}
\newtheorem{corollary}{Corollary}[section]
\def\ep{\varepsilon}
\def\Rn{\mathbb{R}^{n}}
\def\Rm{\mathbb{R}^{m}}
\def\E{\mathbb{E}}
\def\hte{\hat\theta}
\renewcommand{\theequation}{\thesection.\arabic{equation}}
\begin{frontmatter}



\title{Bifurcation and chaotic behaviour in stochastic Rosenzweig-MacArthur prey-predator model
 with non-Gaussian stable L\'evy noise}

\author{Shenglan Yuan\fnref{addr1}}\ead{shenglan.yuan@math.uni-augsburg.de}
\author{Zibo Wang\corref{cor1}\fnref{addr2}}
\ead{zibowang@hust.edu.cn}\cortext[cor1]{Corresponding author}

\address[addr1]{\rm Institut f$\rm\ddot{u}$r Mathematik, Universit$\rm\ddot{a}$t Augsburg,
86135, Augsburg, Germany }
\address[addr2]{\rm Center for Mathematical Sciences, Huazhong
University of Science and Technology, 430074, Wuhan, China}

\begin{abstract}
We perform dynamical analysis on a stochastic Rosenzweig-MacArthur model driven by $\alpha$-stable L\'evy motion. We analyze the existence of the equilibrium points, and provide a clear illustration of their stability. It is shown that the nonlinear model has at most three equilibrium points. If the coexistence equilibrium exists, it is asymptotically stable attracting all nearby trajectories. The phase portraits are drawn to gain useful insights into the dynamical underpinnings of prey-predator interaction. Specifically, we present a transcritical bifurcation curve at which system bifurcates. The stationary probability density is characterized by the non-local Fokker-Planck equation and confirmed by some numerical simulations. By applying Monte Carlo method and using statistical data, we plot a substantial number of simulated trajectories for stochastic system as parameter varies. For initial conditions that are arbitrarily close to the origin, parameter changes in noise terms can lead to significantly different future paths or trajectories with variations, which reflect chaotic behaviour in mutualistically interacting two-species prey-predator system subject to stochastic influence.
\end{abstract}

\begin{keyword}
stochastic Rosenzweig-MacArthur model, prey-predator interaction, transcritical bifurcation, stationary probability density, chaotic behaviour.


\emph{2020 Mathematics Subject Classification}:
37H20, 65P20, 70K50.

\end{keyword}

\end{frontmatter}


\section{Introduction}
The Lotka-Volterra equations \cite{L20,L02,V}, also known as the prey-predator equations, are a pair of first-order nonlinear differential equations, frequently used to describe the dynamics of biological systems in which two species interact, one as prey and the other as predator. The Lotka-Volterra model as the most classical prey-predator model supposes an unlimited food supply for interacting species, while most of interactions occur in limited resource environment. It is more realistic to assume that the interactions would saturate because of the limiting carrying capacity of the environment. The Rosenzweig-MacArthur models with limited resources attract more and more scholars to study the effect of various factors on prey-predator interactions; see various studies \cite{DLM,JAM,SS}. This kind of models can be applied to species living in the world's oceans and animal populations on land \cite{BSD,FRFL,ZKS}. For example, sea lions and penguins, red and grey squirrels, and ants and termites are all species which fall into this category \cite{L}.

Excessive human activities seriously cause climate change and global warming, such as rising temperatures, melting glaciers and sea ice, setting wildlife populations and habitats on the move, and increasing extreme weather events. Now the world is totally different, the surface of this planet is utterly transformed, the extinction speed of the creatures is beyond our imagination. We can find the plastic everywhere, even in the seabird's stomach. Global greenhouse gas emissions are likely to rise to record levels.  In species interactions, the prey hopes to evolve to avoid being caught by the predator, whereas the predator hopes to be able to catch the prey as efficiently as possible. They are inevitably influenced by environmental effects: pollution, refuge, severe drought, overuse of pesticides, drinking water shortages, catastrophic flood, unprecedented burning and other external factors \cite{BS,GHJS,MMIA,PSKD}.

Stochastic noises can mimic the fluctuations in the environment of the dynamical systems \cite{LA,Du15}. The most commonly used stochastic driving process is Gaussian white noise \cite{RP}, but it only describes some fluctuations around mean
value without jumps. Non-Gaussian noise is more close to the reality, which has infinite variance and simulates small perturbations combined with discontinuously unpredictable jumps \cite{Sato}, such as $\alpha$-stable noise \cite{ST}.

The goal of this work is to study the dynamics of a stochastic model as an extension of Rosenzweig-MacArthur model with Holling type III functional response. To the best of our knowledge, the chaotic dynamics of stochastic system \eqref{SRM} have not been studied. Long-term prediction is a challenging yet important task. A description of individual trajectories for stochastic system \eqref{SRM} is not so good, but a statistical description is more appropriate. When the noise intensities are large, stochastic perturbations strong enough to produce a pronounced effect on the dynamical behavior of the model \eqref{SRM} and induce chaos. With the stability indexes decreasing, there are more and more big jumps of $\alpha$-stable noises having the potential to cause abrupt changes. Hence, the trajectories may become chaotic.

We outline the format of this paper as follows. In Section \ref{Sb}, we determine that the model \eqref{RM} has three possible equilibrium points including the conditions for their existence and stability properties. The trivial equilibrium point $Z_1$ is always unstable while
two other equilibrium points, i.e., the predator extinction point $Z_2$ and the coexistence point $Z_3$, are
conditionally stable. We numerically demonstrate the stabilities of the equilibrium points, and carefully consider the occurrence of transcritical bifurcation. In Section \ref{Ss}, we establish the non-local Fokker-Planck equation for stochastic Rosenzweig-MacArthur model \eqref{SRM},  whose solution is stationary density function exhibited by stereoscopic graphs. In Section \ref{Cd}, we discuss chaotic dynamics of stochastic system \eqref{SRM} using solution curves and phase-space diagrams. Several numerical simulations are also given to graphically display the dynamical complexities and pattern of the populations in this system. We end our work with a brief conclusion including important stepping stones to future research in Section \ref{Cf}.

\section{Rosenzweig-MacArthur's model }\label{Sb}
The Rosenzweig-MacArthur's prey-predator system \cite{RM} builds upon the Lotka-Volterra model, adding realism with both logistic growth of the prey, and a limit to the consumption rate of the predator,
\begin{eqnarray}\label{RM}
\left\{\begin{array}{l}
 \dot{X}=rX-cX^2-Yf(X),  \\
\dot{Y}=(Ef(X)-\mu)Y,
 \end{array}\right.
\end{eqnarray}
where $X\geq0$ represents the number of the prey population, and $Y\geq0$ is the size of the predator population. The Holling type III functional response
\begin{equation*}
f(X)=\frac{CX^2}{X^2+k^2}
\end{equation*}
describes a nonlinear consumption, which grows with respect to $X$ when $X$ is small, saturates at the maximum food intake $C$ when $X$ is large.
Model \eqref{RM} represents an interaction between two populations with a prey-predator relationship.
The ecological parameters are all positive constants as described in Table \ref{P}.
\begin{table}[H]
  \centering
\begin{tabular}{|c|c|}
\hline
\text{Parameters}& \text{Description}\\
\hline
$c=0.02$ &\text{the competition factor of prey}\\
\hline
$E=0.4$ &\text{the assimilation efficiency of predator}\\
\hline
$C=1$ &\text{the maximum food intake of predator}\\
\hline
 $k=10$ &\text{the half-saturation constant of functional response}\\
\hline
$r\in(0,2.5)$ &\text{the intrinsic growth rate of prey}\\
\hline
$\mu\in(0,0.22)$ &\text{the mortality rate of predator}\\
\hline
$\sigma_i\in(0,1)$&\text{the intensities of noise}\\
\hline
$\alpha_i\in(0,2)$&\text{the indexes of stability}\\
\hline
\end{tabular}
  \caption{Parameters of the system \eqref{SRM} }\label{P}
\end{table}
\subsection{Stationary states}
The system \eqref{RM} has at most three equilibria by using the equilibrium equations
\begin{align*}
 rX-0.02X^2-Yf(X)&=0, \\
  (0.4f(X)-\mu)Y&=0.
\end{align*}
A trivial zero population solution $Z_{1}=(0,0)$ and a prey-only solution $Z_{2}=(50r,0)$ always exist for all parameter settings. The local stability of all equilibrium points can be studied from the linearization of system \eqref{RM}.
Linearize by calculating the Jacobian matrix
\begin{equation*}
 J=\left(\begin{array}{cc}
     r-0.04X-Yf'(X) & -f(X) \\
     0.4Yf'(X) & 0.4f(X)-\mu
   \end{array}\right),
\end{equation*}
where $f'(X)=C\frac{200X}{(X^{2}+100)^2}$. Linearize at the origin equilibrium $Z_1=(0,0)$:
\begin{equation*}
 J_{Z_1}=\left(\begin{array}{cc}
     r & 0 \\
     0 & -\mu
   \end{array}\right).
\end{equation*}
The eigenvalues of $J_{Z_1}$ are $\lambda_1=r>0$ and $\lambda_2=-\mu<0$. Therefore, the trivial equilibrium is an unstable saddle point.
The Jacobian matrix for the predator extinction equilibrium $Z_2=(50r,0)$ is:
\begin{equation*}
 J_{Z_2}=\left(\begin{array}{cc}
     -r & -f(50r) \\
     0 & 0.4f(50r)-\mu
   \end{array}\right).
\end{equation*}
The prey-only equilibrium has eigenvalues $\lambda_1=-r<0$ and $\lambda_2=0.4f(50r)-\mu$.
Hence, it is a stable node (locally asymptotically stable) based on that all two of the eigenvalues are
negative when $\lambda_2<0$, i.e.,
\begin{equation}\label{Crm}
\frac{0.4r^{2}}{0.04+r^{2}}<\mu,
\end{equation}
is saddle point when $\lambda_2>0$, and undergoes a transcritical
bifurcation whenever $\lambda_2=0$. The appearance of transcritical bifurcation is caused by the changing of the sign of $\lambda_2$.
The other positive equilibrium in the first quadrant corresponds to a stationary coexistence of prey and predator,
and satisfies the following conditions:
\begin{align}\label{H}
  Y&=\frac{rX-0.02X^2}{f(X)}=\frac{(r-0.02X)(X^{2}+100)}{X}, \\ \label{hP}
  g(X)&=0.4f(X)-\mu=\frac{0.4X^{2}}{X^{2}+100}-\mu=0.
\end{align}
Note that condition \eqref{hP} gives the $X-$component of the coexistence equilibrium solution.
What is more, the condition \eqref{hP} and thus the $X-$component are $r-$independent. The
net per-capita predator growth $g(X)$ equals $-\mu$ for $X=0$, which is strictly increasing and levels off at $0.4-\mu>0$ for large $X$.
Thus, Eq. \eqref{hP} has only one positive root, i.e., the $X-$component of the coexistence 
equilibrium. Rewriting Eq. \eqref{hP} into $(0.4-\mu)X^{2}=100\mu$, which yields
\begin{equation*}
X=10\sqrt{\frac{\mu}{0.4-\mu}}.
\end{equation*}
Substituting this into Eq. \eqref{H} gives
\begin{equation*}
Y=\frac{4r\sqrt{0.4-\mu}-0.8\sqrt{\mu}}{(0.4-\mu)\sqrt{\mu}}.
\end{equation*}
The $Y-$component is positive if and only if
\begin{equation}\label{C}
\frac{0.4r^{2}}{0.04+r^{2}}>\mu,
\end{equation}
which indicates the all of species coexist since the $X-$component is already positive.  If the condition \eqref{C} is fulfilled, then all three equilibrium points of system \eqref{RM} exist. We remark that when \eqref{C} does not hold, the coexistence equilibrium point $Z_3$ does not exist in this case.

We explicitly express and analyze the coexistence equilibrium point
\begin{equation*}
Z_3=(X,Y)=\Big(10\sqrt{\frac{\mu}{0.4-\mu}},\frac{4r\sqrt{0.4-\mu}-0.8\sqrt{\mu}}{(0.4-\mu)\sqrt{\mu}}\Big).
\end{equation*}
The Jacobian matrix evaluated at the prey-predator equilibrium $Z_3$ is
\begin{equation*}
 J_{Z_3}=\left(\begin{array}{cc}
     r-0.4\sqrt{\frac{\mu}{0.4-\mu}}-5r(0.4-\mu)+\sqrt{\mu(0.4-\mu)} & -\frac{\mu}{0.4} \\
     2r(0.4-\mu)-0.4\sqrt{\mu(0.4-\mu)} & 0
   \end{array}\right).
\end{equation*}
The eigenvalues of the Jacobian matrix $J_{Z_3}$ are the solutions of the characteristic equation
\begin{equation*}
\lambda^{2}-\Big(r-0.4\sqrt{\frac{\mu}{0.4-\mu}}-5r(0.4-\mu)+\sqrt{\mu(0.4-\mu)}\Big)\lambda+\mu\big(5r(0.4-\mu)-\sqrt{\mu(0.4-\mu)}\big)=0.
\end{equation*}
Solving this quadratic equation for $\lambda$ obtains the roots as
\begin{align*}
\lambda_{1,2}&=\frac{1}{2}\Big(r-0.4\sqrt{\frac{\mu}{0.4-\mu}}-5r(0.4-\mu)+\sqrt{\mu(0.4-\mu)}\Big) \\
             &~~~~\pm\sqrt{\frac{1}{4}\Big(r-0.4\sqrt{\frac{\mu}{0.4-\mu}}-5r(0.4-\mu)+\sqrt{\mu(0.4-\mu)}\Big)^{2}-\mu\big(5r(0.4-\mu)-\sqrt{\mu(0.4-\mu)}\big)}\\
             &:=\varphi\pm i\psi.
\end{align*}
Since $\varphi<0$ and $\psi>0$ for parameters $r\in(0,2.5)$ and $\mu\in(0,0.22)$ in Table \ref{P}, the positive coexistence equilibrium point $Z_3$ is a stable spiral attracting all closer enough trajectories.
\begin{figure}[H]
\begin{center}
  \begin{minipage}{2.1in}
\leftline{(a)}
\includegraphics[width=2.1in]{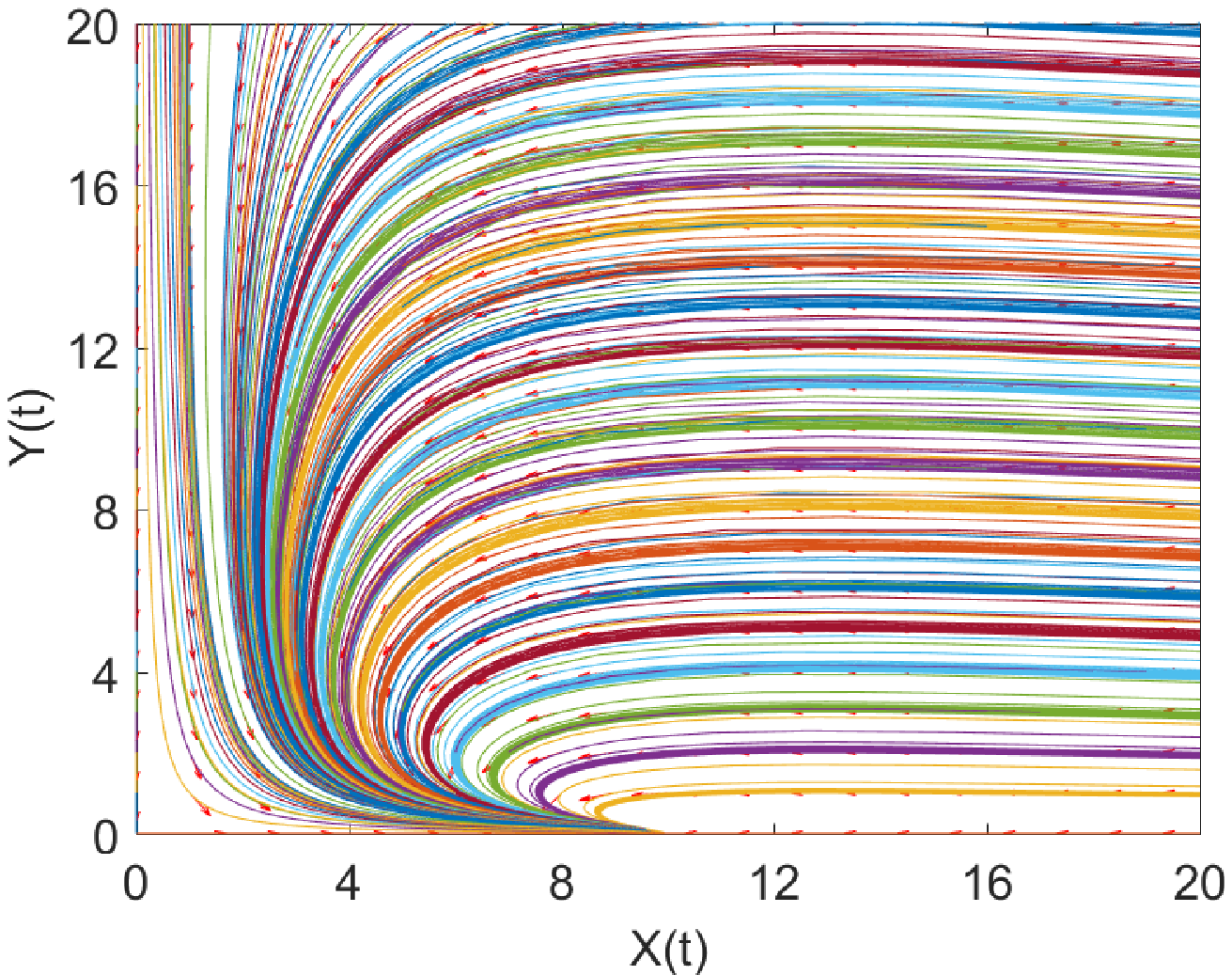}
\end{minipage}
\hfill
\begin{minipage}{2.1in}
\leftline{(b)}
\includegraphics[width=2.1in]{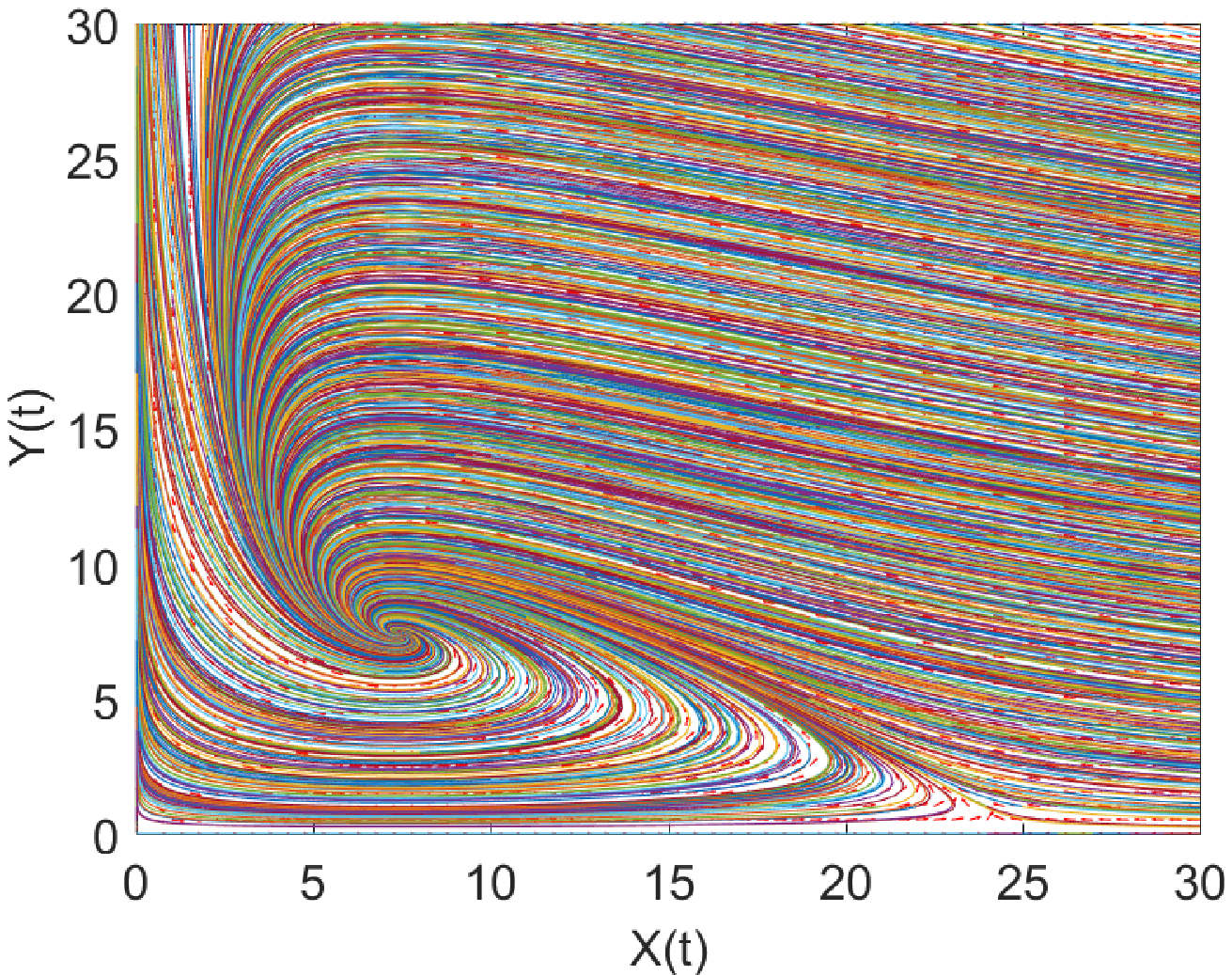}
\end{minipage}
\hfill
  \begin{minipage}{2.1in}
\leftline{(c)}
\includegraphics[width=2.1in]{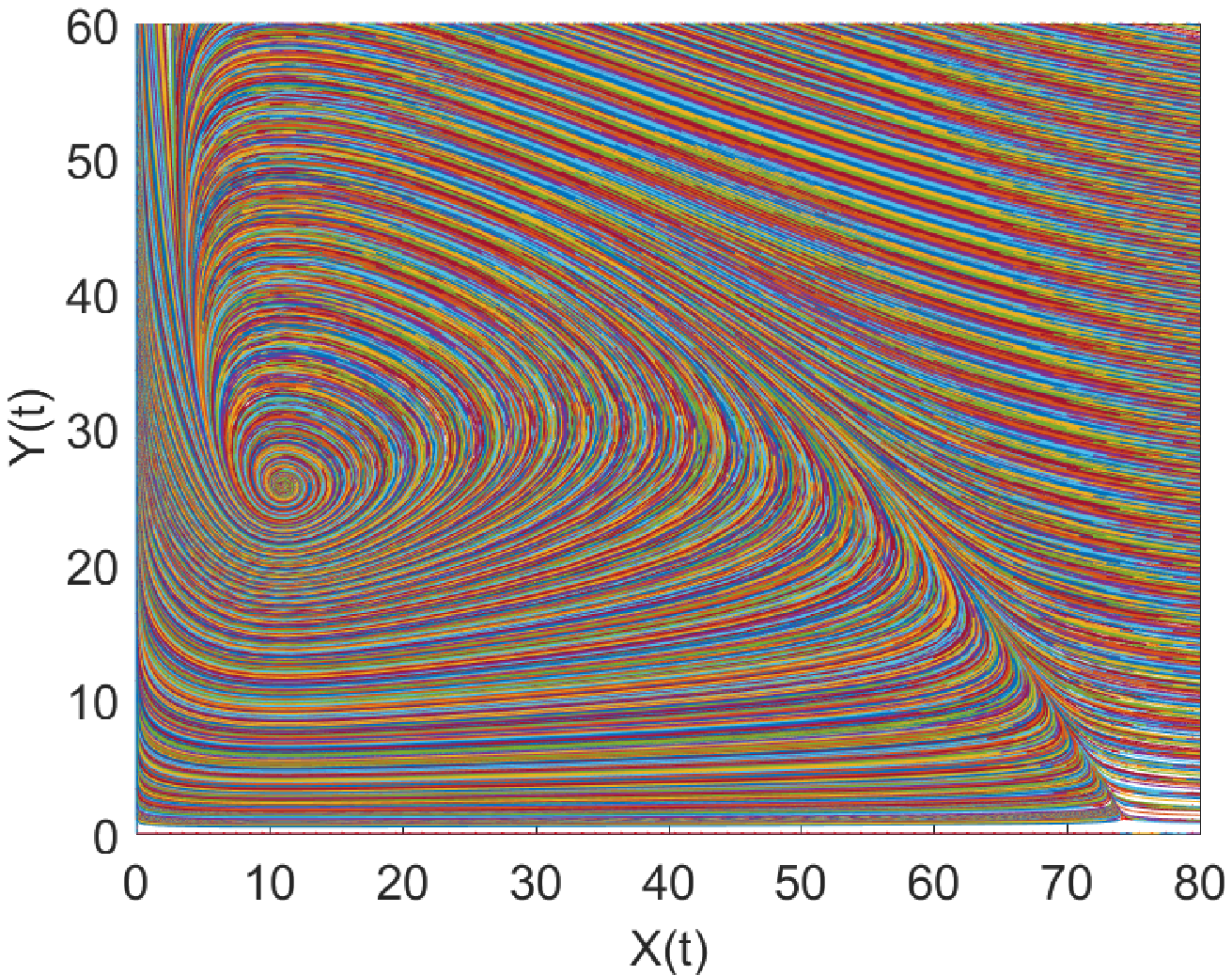}
\end{minipage}
\caption{Phase portraits showing various trajectories of model \eqref{RM} under different parameters: (a) $r=0.2$, $\mu=0.25$; (b) $r=0.5$, $\mu=0.14$; (c) $r=1.5$, $\mu=0.22$.}\label{Fig1}
\end{center}
\end{figure}

Our analytical findings in this subsection are justified by performing numerical simulations. In Fig \ref{Fig1}(a), we plot a phase plane portrait. The parameter values $r=0.2$ and $\mu=0.25$ lead to the non-existence of equilibrium point $Z_3$, i.e., there is no coexistence equilibrium. The predator extinction point  $Z_2=(10,0)$ is an asymptotically stable equilibrium point since the eigenvalues $\lambda_1=-0.2$ $\lambda_2=-0.05$ are both negative. It is seen that all solutions with different initial values are convergent to the stable node $Z_2$. For the extinction equilibrium $Z_1=(0,0)$, there is no population.

The phase plane diagram in Fig \ref{Fig1}(b) shows that the coexistence point $Z_3$ is asymptotically stable when $r=0.5$ and $\mu=0.14$. The eigenvalues of $J_{Z_3}$ are given by $\lambda_{1,2}\approx-0.13\pm i0.22$. Therefore, the fixed point at $Z_3=(7.34,7.4)$ is a stable node or spiral and all phase paths inside the first quadrant end up in $Z_3$. Besides $Z_1=(0,0)$, the predator-free equilibrium $Z_2=(25,0)$ is also a saddle point because of $\lambda_2=0.2>0$. In Fig \ref{Fig1}(c), we show the phase portrait for the case of $r=1.5$ and $\mu=0.22$. We see that the boundary equilibria $Z_1=(0,0)$ and $Z_2=(75,0)$ are saddle points, while the unique interior equilibrium $Z_3=(11.06,25.7)$ is asymptotically stable.
All trajectories lying in the first quadrant are drawn to the fixed point $Z_3$ no matter what the initial values of $X(t)$ and $Y(t)$.

\subsection{Bifurcation analysis}
We have a detailed discussion on the condition \eqref{C} for the existence of the prey-predator equilibrium point $Z_3$.
It is noted that the stability condition \eqref{Crm} of $Z_2$ contradicts the existence condition \eqref{C} for the coexistence point $Z_3$. Consequently, if $Z_2$ is asymptotically stable, then the interior equilibrium $Z_3$ does not exist. Those conditions also indicate the existence of transcritical bifurcation.

Based on the existence and stability results of equilibrium points of the system \eqref{RM},  substitution of $X=50r$ into Eq. \eqref{hP} provides a curve of transcritical bifurcation in the $(r,\mu)$ parameter plane
\begin{equation*}
T=\Big\{(r,\mu): \mu=\frac{0.4r^{2}}{0.04+r^{2}}\Big\}.
\end{equation*}
When $(r,\mu)$ passes through $T$, the equilibrium point $Z_2$ undergoes a transcritical bifurcation, which changes from a sink to a source as one eigenvalue of the Jacobian matrix $J_{Z_2}$ changes sign from negative to positive.

Now we discuss the existence of the interior attractor $Z_3$ by considering
the regions divided by $T$. In addition to the trivial equilibrium $Z_1$ and the unstable predator-free equilibrium $Z_2$,  system \eqref{RM} also has the asymptotically stable coexistence point $Z_3$ if $(r,\mu)$ in the region below $T$. The predator extinction point $Z_2$ is asymptotically stable if $(r,\mu)$ in the region up $T$. This means that prey will survive in the system \eqref{RM}, while predator will go extinct.

\begin{figure}[H]
\begin{center}
\begin{minipage}{2.1in}
\leftline{(a)}
\includegraphics[width=2.1in]{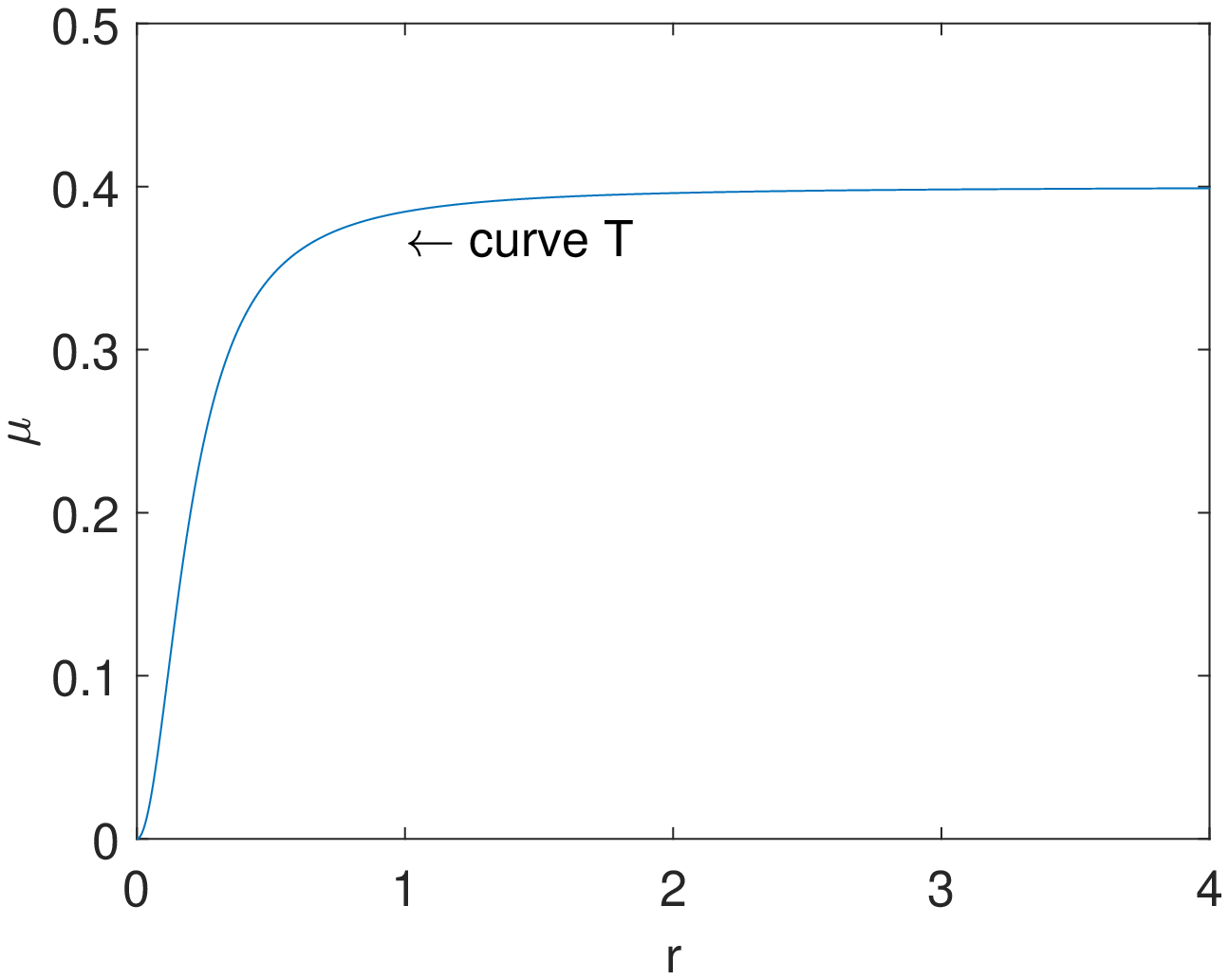}
\end{minipage}
\hfill
  \begin{minipage}{2.1in}
\leftline{(b)}
\includegraphics[width=2.1in]{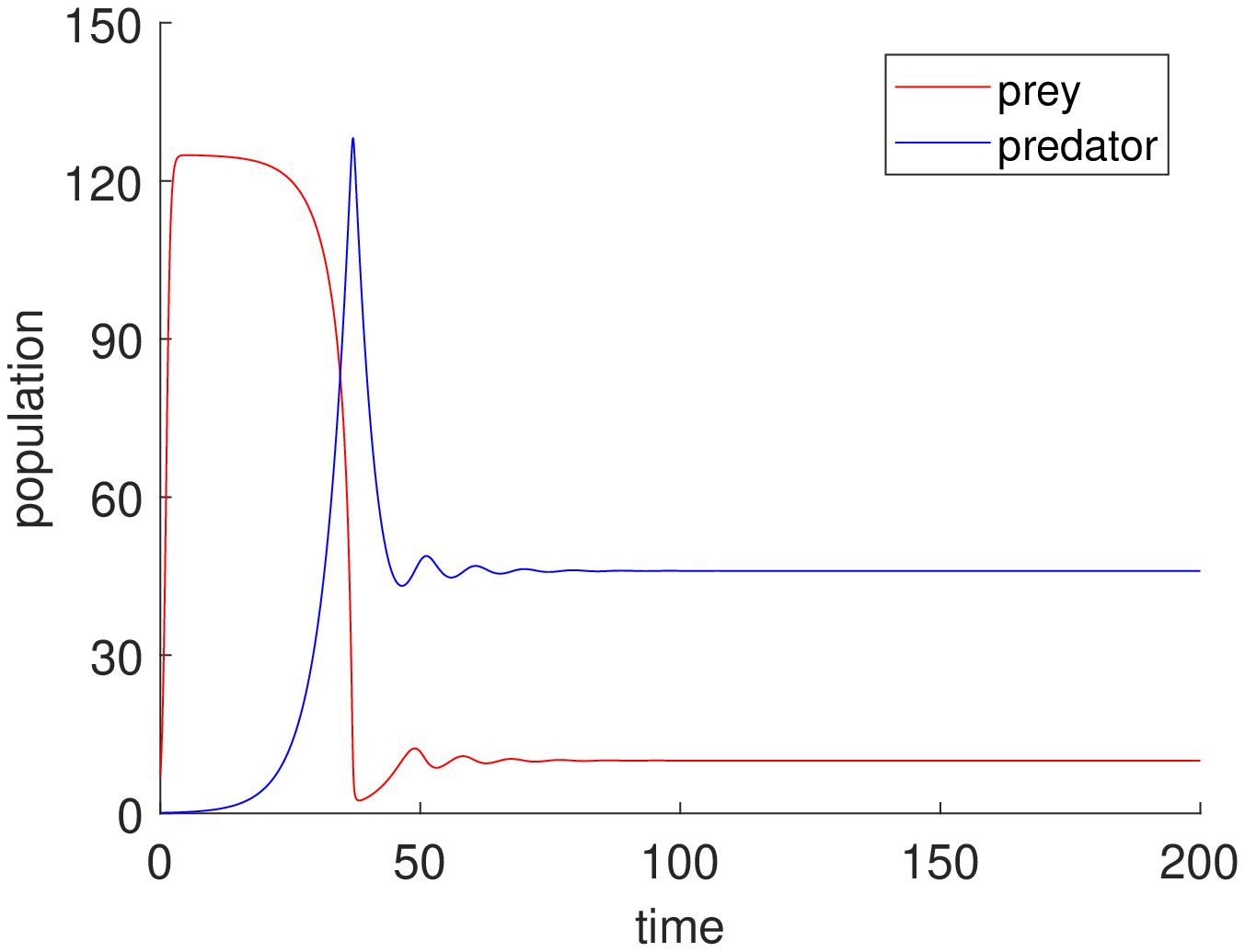}
\end{minipage}
\hfill
  \begin{minipage}{2.1in}
\leftline{(c)}
\includegraphics[width=2.1in]{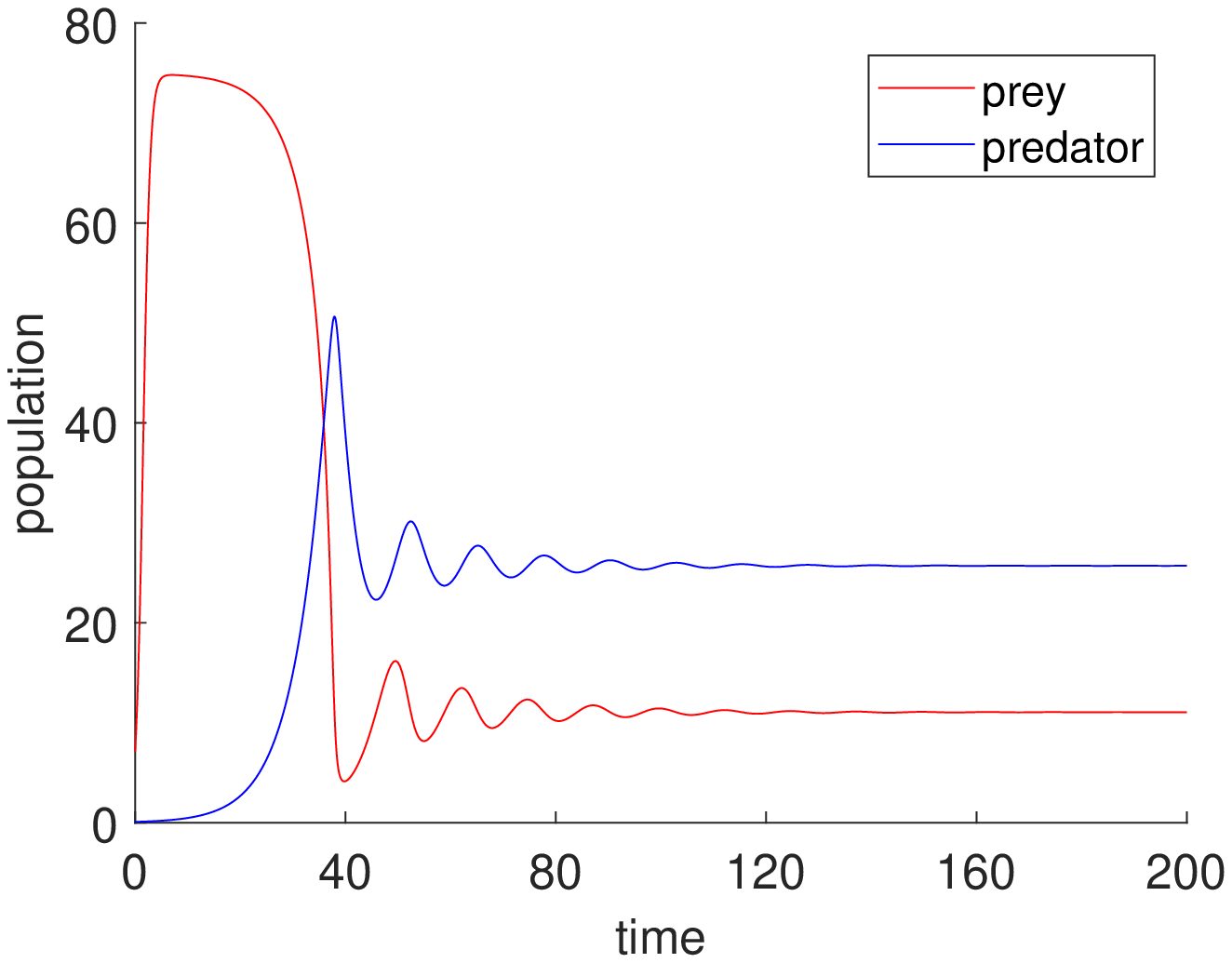}
\end{minipage}
\caption{(a) Curve $T$ emerges from the origin and levels off at $\mu=0.4$ from large $r$. The prey-predator equilibrium $Z_3$ that bifurcates from the prey-only equilibrium $Z_2$ exists below $T$; (b) As time evolves, the populations $X(t)$ and $Y(t)$ stabilize at $Z_3=(10,46)$ with parameters $r=2.5$ and $\mu=0.2$; (c) For $r=1.5$ and $\mu=0.22$, the prey $X(t)$ and predator $Y(t)$ are oscillating in the beginning stage, but move toward $11.06$ and 25.7 respectively.}\label{Fig2}
\end{center}
\end{figure}

The system \eqref{RM} undergoes a transcritical bifurcation at $T$ with stability-instability switching of $Z_2$ and
creation/destruction of $Z_3$. From  Fig \ref{Fig2}(a), we see that the curve $T$ controlled by $r$ and $\mu$ is a threshold: if $\mu$ is greater than $0.4r^{2}/(0.04+r^{2})$, then $Z_2$ is the unique asymptotically stable equilibrium point; if $\mu$ is less than $0.4r^{2}/(0.04+r^{2})$, then $Z_2$ becomes unstable and there appears an asymptotically stable equilibrium point $Z_3$.

To see the detail dynamics of the system \eqref{RM}, we plot the populations of interacting species.  In Fig \ref{Fig2}(b),
the populations rise and fall, eventually settle down to constant values.
Because the condition
\begin{equation*}
\frac{0.4r^{2}}{0.04+r^{2}}\approx0.4>0.2=\mu,
\end{equation*}
is valid, the system is asymptotically convergent to the point $Z_3=(10,46)$ in terms of the parameters $r=2.5$ and $\mu=0.2$.
It is clearly seen in Fig \ref{Fig2}(c) that for $r=1.5$ and $\mu=0.22$, the trajectories of $X(t)$ and $Y(t)$ are
drawn to $11.06$ and 25.7 respectively, and once there, remain there in the long term.
\section{Stochastic system}\label{Ss}
Increasingly, for many application areas, it is becoming important to include elements of nonlinearity and non-Gaussianity in order to model accurately the underlying dynamics of a dynamical system by using stochastic differential equation modelling techniques \cite{YZD}. The dynamics of Rosenzweig-MacArthur model perturbed by $\alpha$-stable L\'evy noise can be represented mathematically with two nonlinear stochastic differential equations in $\mathbb{R}$, given by
\begin{eqnarray}\label{SRM}
\left\{\begin{array}{l}
 dX=\big(rX-0.02X^2-Yf(X)\big)dt+\sigma_1dL_{t}^{\alpha_1},  \\
 dY=\big[(0.4f(X)-\mu)Y\big]dt+\sigma_2dL_{t}^{\alpha_2}.
 \end{array}\right.
\end{eqnarray}
 For the prey and predator $X, Y>0$, the two populations oscillate. Both populations are influenced by external fluctuations. The stochastic noise terms $\{L_t^{\alpha_{i}}: t\geq0\}, i=1,2$ are independent real-valued non-Gaussian symmetric $\alpha$-stable processes with L\'evy triplets $(0,0,\nu_{\alpha_{i}})$ on probability spaces $(\Omega^{i},\mathcal{F}^{i},\mathbb{P}^{i})$. Note that $L_{t}^{\alpha}=\big(L_{t}^{\alpha_1}, L_{t}^{\alpha_2}\big)^{\top}$ is a two-dimensional $\alpha$-stable L\'evy process with the L\'evy triplet $(l,Q,\nu_\alpha)$, where $l=\big(0, 0\big)^{\top}$, $Q$ is $2\times2$ null matrix, and $\nu_\alpha(du,dv)=\nu_{\alpha_{1}}(du)\delta_{0}(dv)+\nu_{\alpha_{2}}(dv)\delta_{0}(du)$.
The L\'evy measure $\nu_{\alpha_1}$ satisfies $\int_{\mathbb{R}\setminus{\{0\}}}(|u|^{2}\wedge1)\nu_{\alpha_1}(du)<\infty$, which is determined by
\begin{equation*}
\nu_{\alpha_1}(du)=c(1,{\alpha})\frac{1}{|u|^{1+\alpha}}dy~~~~~~\text{and}~~~~~~c(1,\alpha)=\frac{\alpha\Gamma(\frac{1+\alpha}{2})}{{2^{1-\alpha}\sqrt{\pi}}\Gamma{(1-\frac{\alpha}{2})}},
\end{equation*}
 where $\Gamma$ is the Gamma function. The L\'evy measure $\nu_{\alpha_2}$ is similarly defined.

A knowledge of the stationary probability density gives us a wealth of statistical information in the asymptotic regime \cite{DSS,HLS}.
Now we study how an ensemble of initial conditions, characterized by an initial density $p(x_0,y_0,0)$, propagates under the action of stochastic system \eqref{SRM}. The evolution of this density $p(x,y,t)$ is governed by the non-local Fokker-Planck equation:
\begin{align}\nonumber
	\frac{\partial}{\partial t}p(x,y,t)&=\big(0.04x+yf'(x)+\mu-0.4f(x)-r\big)p(x,y,t)  \\ \nonumber
	&~~~~+\big(0.02x^{2}+yf(x)-rx\big)\frac{\partial}{\partial x}p(x,y,t)+y(\mu-0.4f(x))\frac{\partial}{\partial y}p(x,y,t)  \\ \nonumber
	&~~~~+\sigma_1^{\alpha_1}\int_{\mathbb{R}\backslash\{0\}}(p(x+u,y,t)-p(x,y,t))\nu_{\alpha_1}(du)\\ \label{FP}
&~~~~+\sigma_2^{\alpha_2}\int_{\mathbb{R}\backslash\{0\}}(p(x,y+v,t)-p(x,y,t))\nu_{\alpha_2}(dv).
\end{align}
All calculations obtaining the non-local Fokker-Planck equation \eqref{FP} can be found in the
Appendix. This propagation of the probability density function $p(x,y,t)$ is only a conceptual solution of Eq. \eqref{FP},
it cannot be determined analytically. We solve the Fokker-Planck equation for stationary solution at the stochastic steady state numerically since $p(x,y,t)$ embodies all
available statistical information.
\begin{figure}[H]
\begin{center}
  \begin{minipage}{2.1in}
\leftline{(a)}
\includegraphics[width=2.1in]{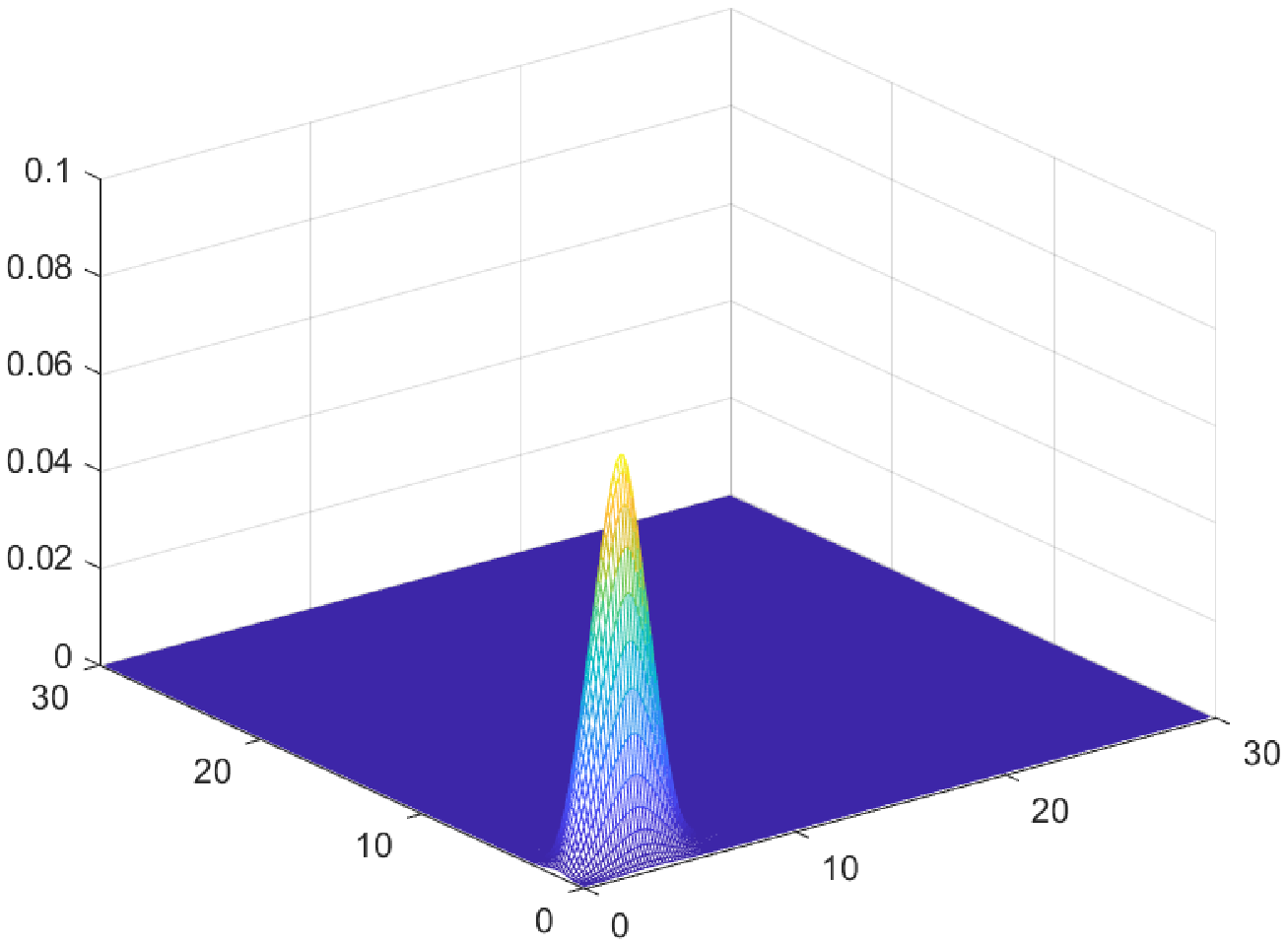}
\end{minipage}
\hfill
\begin{minipage}{2.1in}
\leftline{(b)}
\includegraphics[width=2.1in]{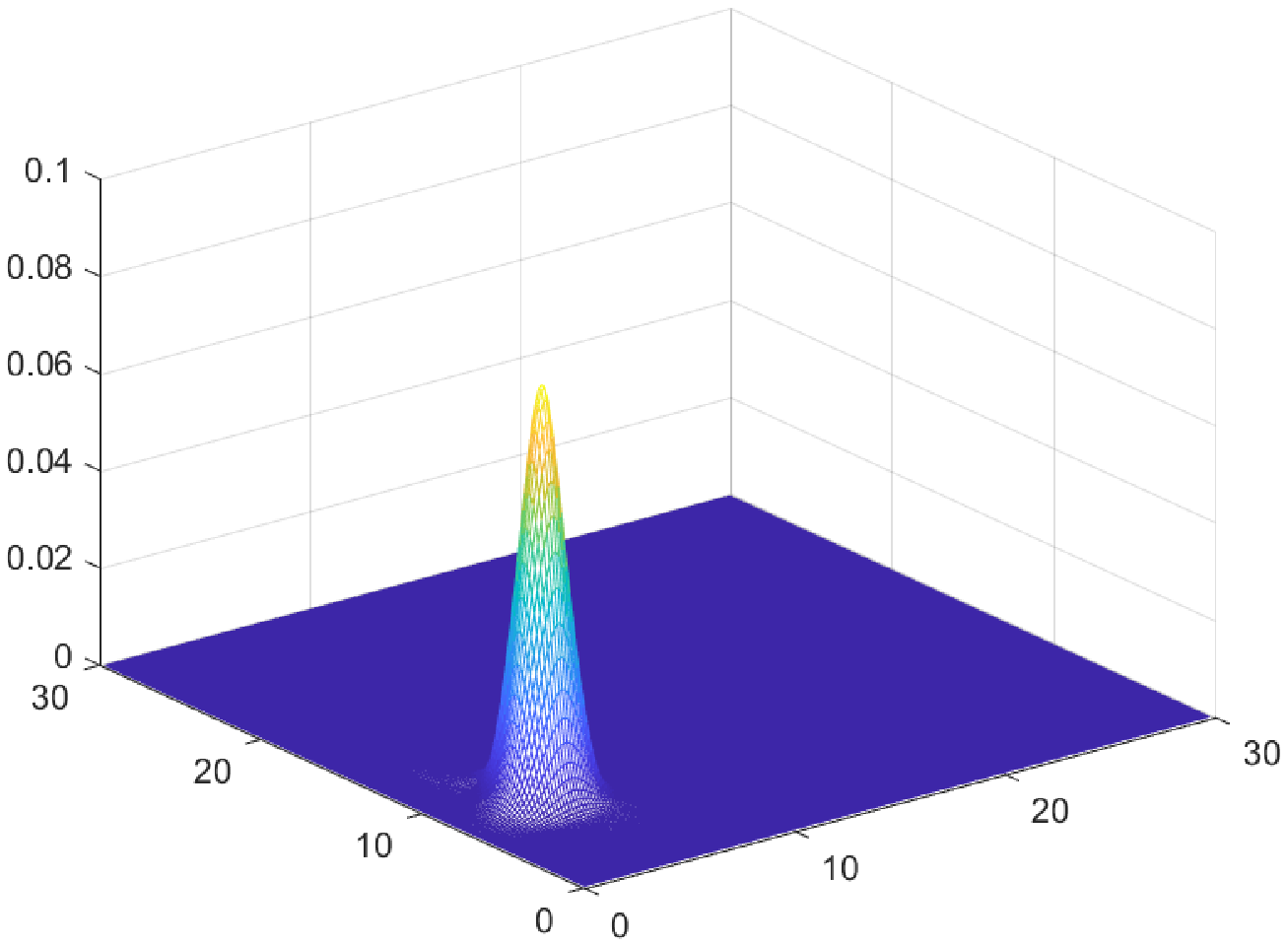}
\end{minipage}
\hfill
  \begin{minipage}{2.1in}
\leftline{(c)}
\includegraphics[width=2.1in]{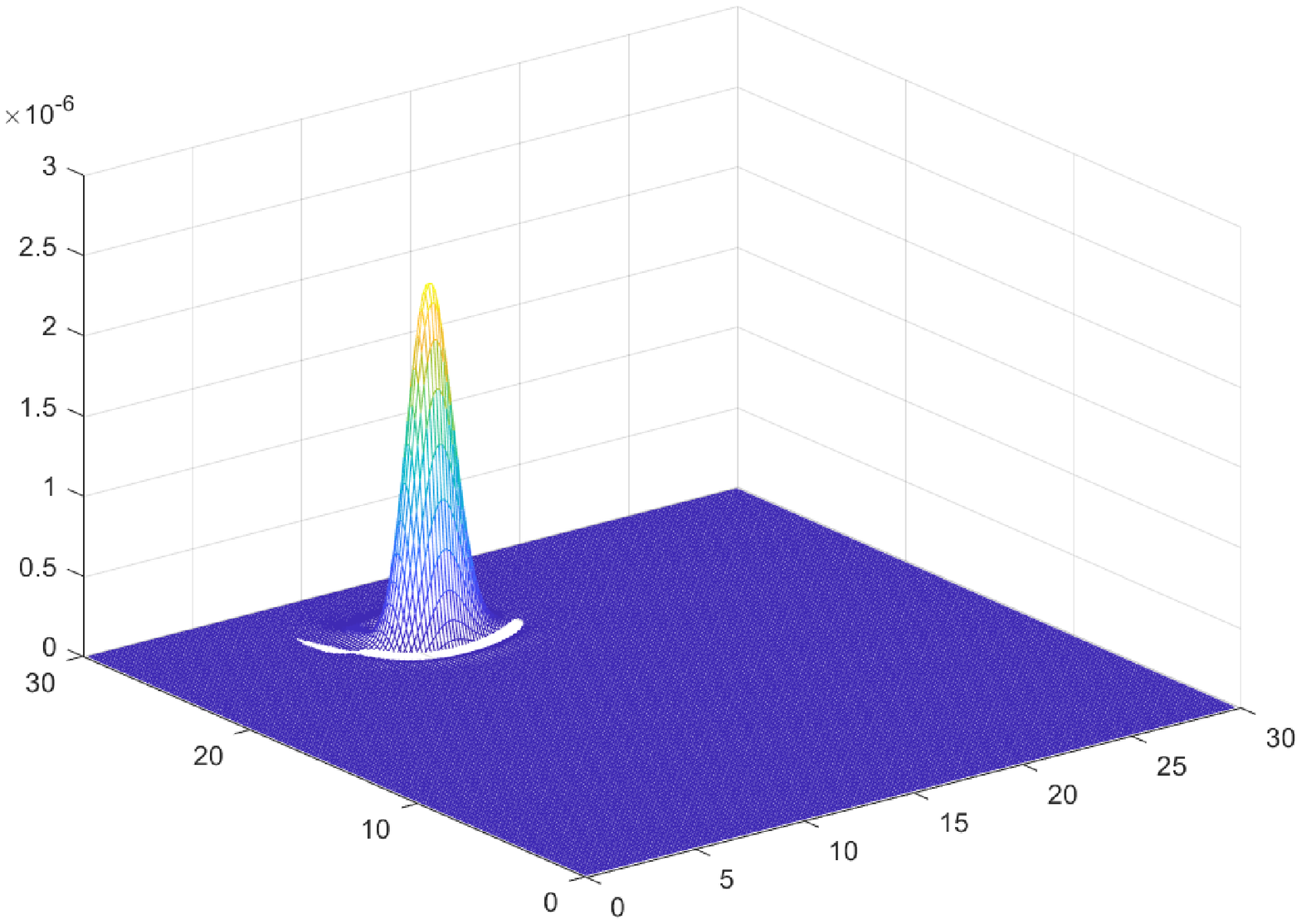}
\end{minipage}
\caption{The probability density functions as stationary solutions of the Fokker-Planck equation \eqref{FP} for noise intensities $\sigma_1=\sigma_2=0.1$ and stability indexes $\alpha_1=\alpha_2=1.5$ with parameters: (a) $r=0.2$, $\mu=0.25$; (b) $r=0.5$, $\mu=0.14$; (c) $r=1.5$, $\mu=0.22$.}\label{Fig4}
\end{center}
\end{figure}

As indicated in Fig \ref{Fig4}, the position of the peak of the probability density function is different because of the changes in the values of
$r$ and $\mu$. The monomodal peak pattern corresponds to the unique steady state as plotted in Fig \ref{Fig1}. The height of the peak of the probability density function for $r=0.2$ and $\mu=0.25$ in Fig \ref{Fig4}(a) is almost the same as
that for $r=0.5$ and $\mu=0.14$ in Fig \ref{Fig4}(b). While the height of the peak of the probability density function for $r=1.5$ and $\mu=0.22$ in Fig \ref{Fig4}(c) is far lower than that for the first two scenarios.

\section{Chaotic dynamics }\label{Cd}
The dynamical properties of the model \eqref{RM}  do not persist if external noises are added to the right-hand sides of the differential equations.
System \eqref{SRM} with $\alpha$-stable L\'evy noises describes the dynamics of the populations as well as
their interactions. The system \eqref{SRM} is not robust (or structurally stable) since small perturbations do affect the qualitative behavior. For both prey and predator subjected to the effects of noises, the paths are more complicated and unpredictable than that of the circumstance where only the prey population (or the predator population) is catalyzed by stochastic noise. To perform a more detailed analysis, we make use of the Monte Carlo simulations to investigate the effects of noise intensities and stability indexes.
For clarity, throughout this section, we will fix the following
parameter quantities: $r=1.5$ and $\mu=0.22$.

\subsection{Effects of noise intensities}
\begin{figure}[H]
\begin{center}
  \begin{minipage}{2.1in}
\leftline{(a)}
\includegraphics[width=2.1in]{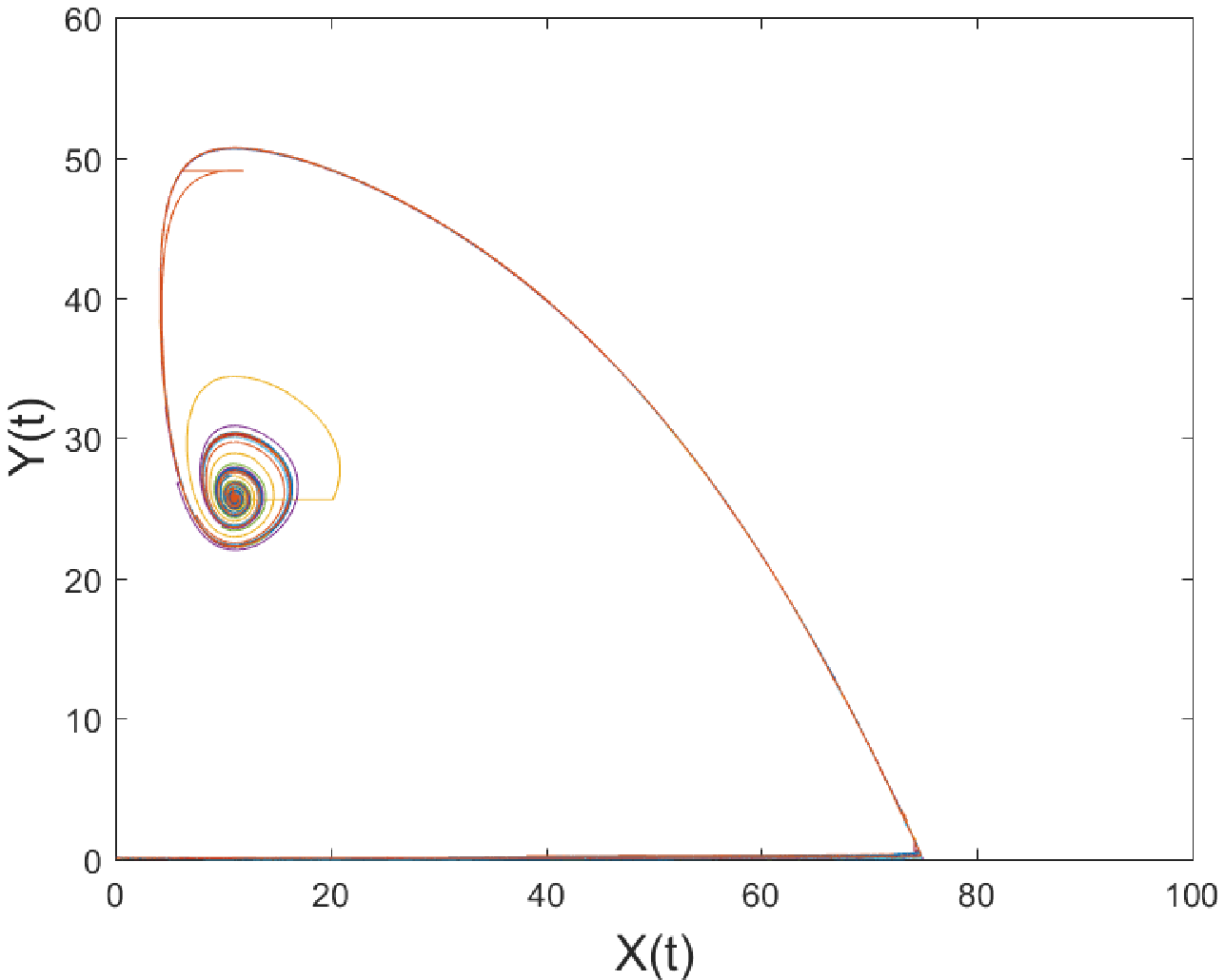}
\end{minipage}
\hfill
\begin{minipage}{2.1in}
\leftline{(b)}
\includegraphics[width=2.1in]{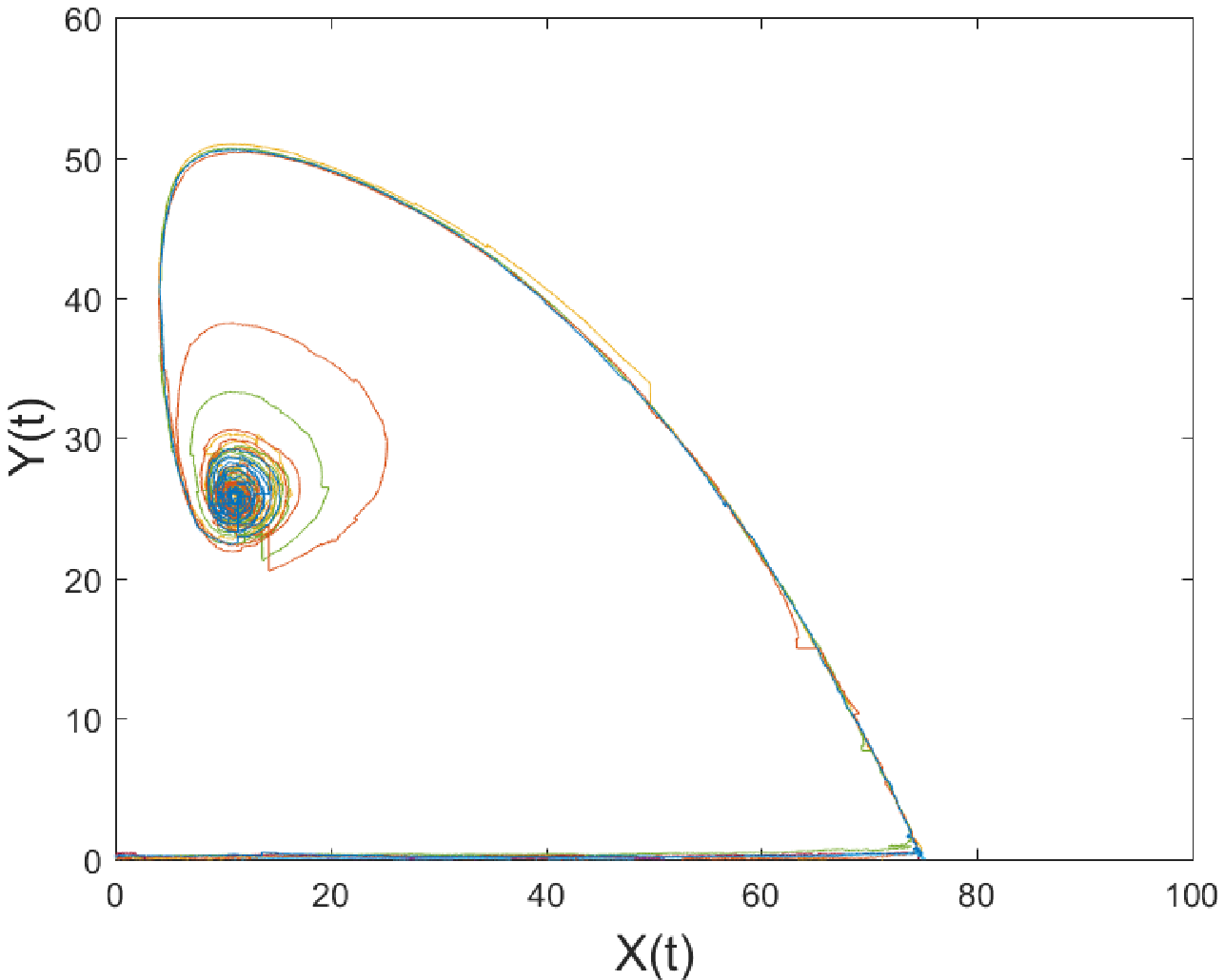}
\end{minipage}
\hfill
  \begin{minipage}{2.1in}
\leftline{(c)}
\includegraphics[width=2.1in]{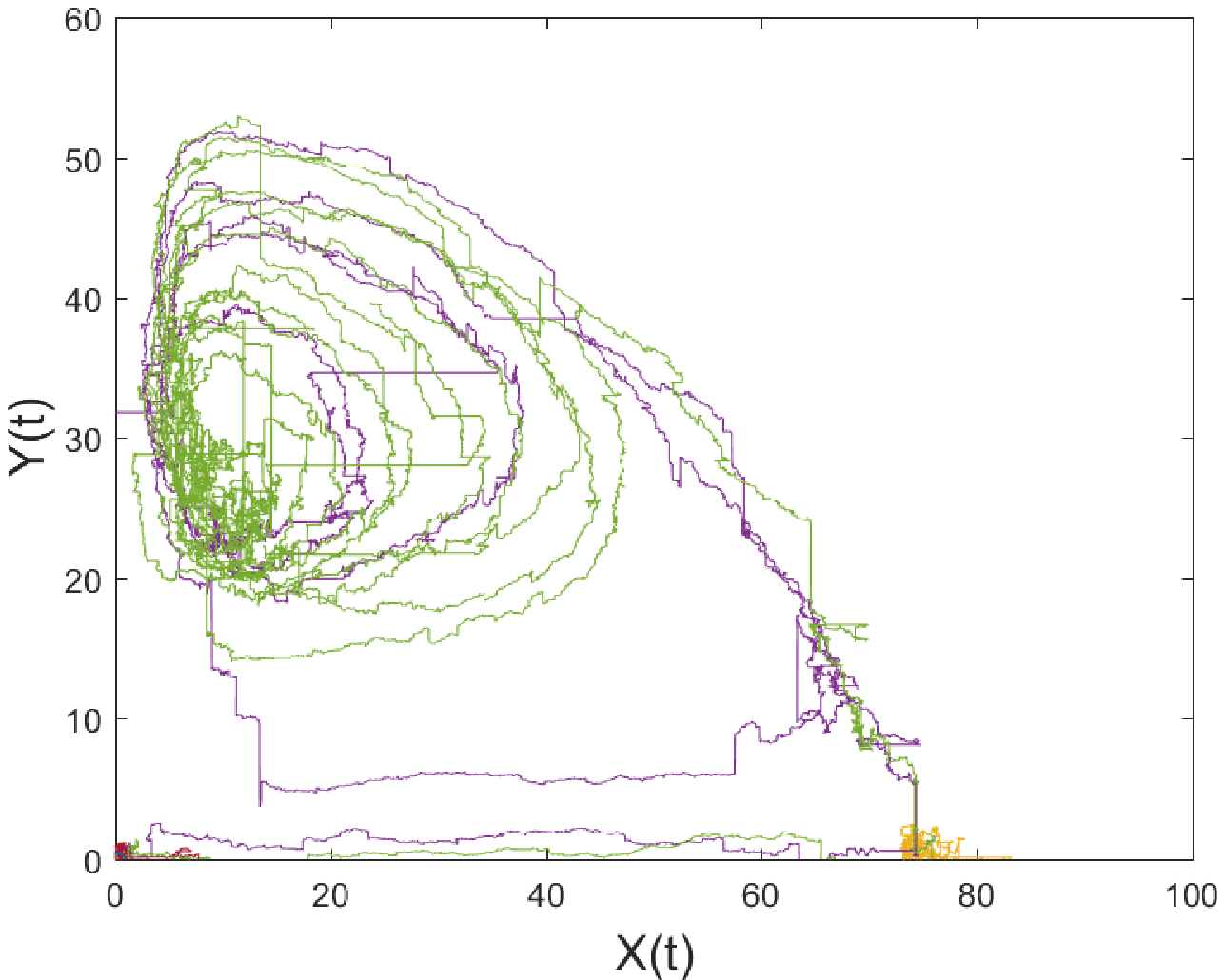}
\end{minipage}
\caption{The interaction of two species at
$\alpha_1=\alpha_2=1.5$: (a) $\sigma_1=\sigma_2=0.01$; (b) $\sigma_1=\sigma_2=0.1$; (c) $\sigma_1=\sigma_2=0.9$.}\label{Fig5}
\end{center}
\end{figure}

 To verify the complexity of population dynamics for system \eqref{SRM} more precisely, we study the interacting species at which the prey and predator populations are subject to unknown disturbances modeled as $\alpha$-stable L\'evy noises. In our computations, we set the stability indexes $\alpha_1=\alpha_2=1.5$. We assume that the noises have equal influence intensities on both the prey $X$ and the predator $Y$. The noises significantly affect the dynamical behaviors of the model \eqref{SRM}. For the small noise intensities $\sigma_1=\sigma_2=0.01$ as in Fig \ref{Fig5}(a), we find that several trajectories describing the interaction of prey-predator converge to the coexistence equilibrium $Z_3=(11.06,25.7)$ which is asymptotically stable. As $\sigma_1=\sigma_2$ are increased toward $0.1$, external noises excite low frequency oscillations of the system paths shown in Fig \ref{Fig5}(b). It can be seen that several winding curves get close to the asymptotically stable equilibrium point $Z_3$. If we strengthen the noise intensities, the curves become sophisticated, confusing and tortuous.  The chaotic behavior for parameter values $\sigma_1=\sigma_2=0.9$ is depicted in Fig \ref{Fig5}(c). The increasing strength of noise intensities can enhance the response of a nonlinear system to external signals.

\begin{figure}[H]
\begin{center}
  \begin{minipage}{2.1in}
\leftline{(a)}
\includegraphics[width=2.1in]{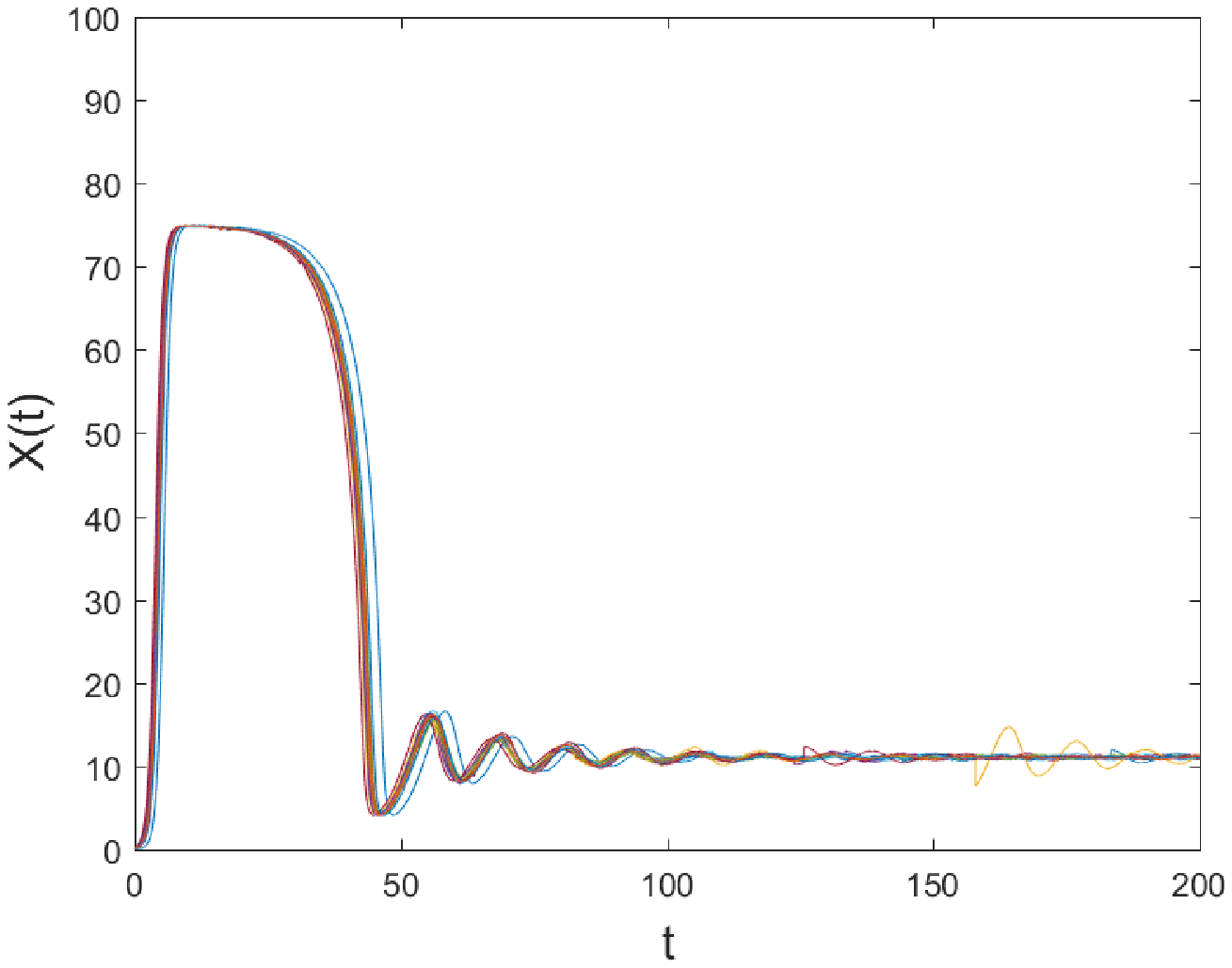}
\end{minipage}
\hfill
\begin{minipage}{2.1in}
\leftline{(b)}
\includegraphics[width=2.1in]{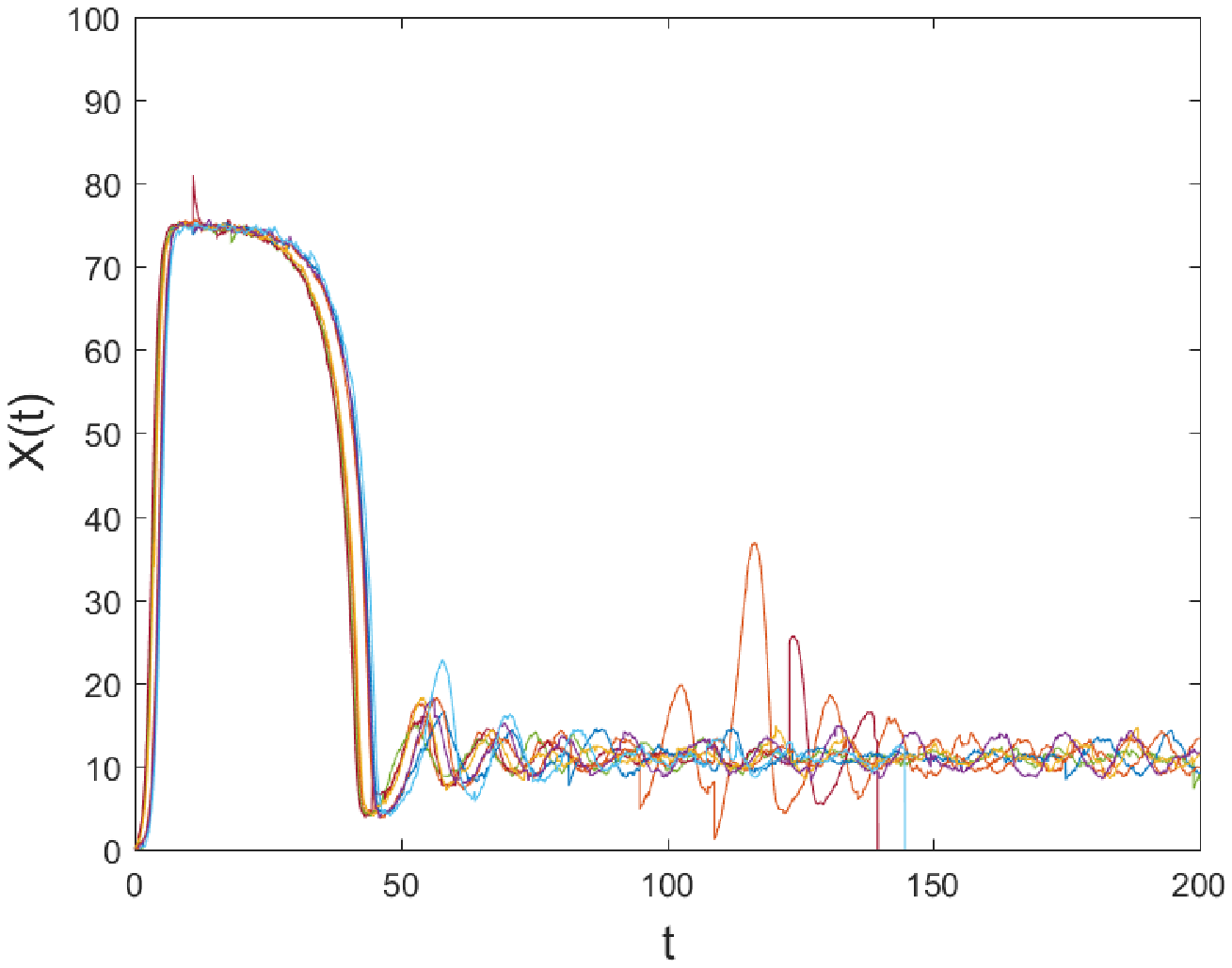}
\end{minipage}
\hfill
  \begin{minipage}{2.1in}
\leftline{(c)}
\includegraphics[width=2.1in]{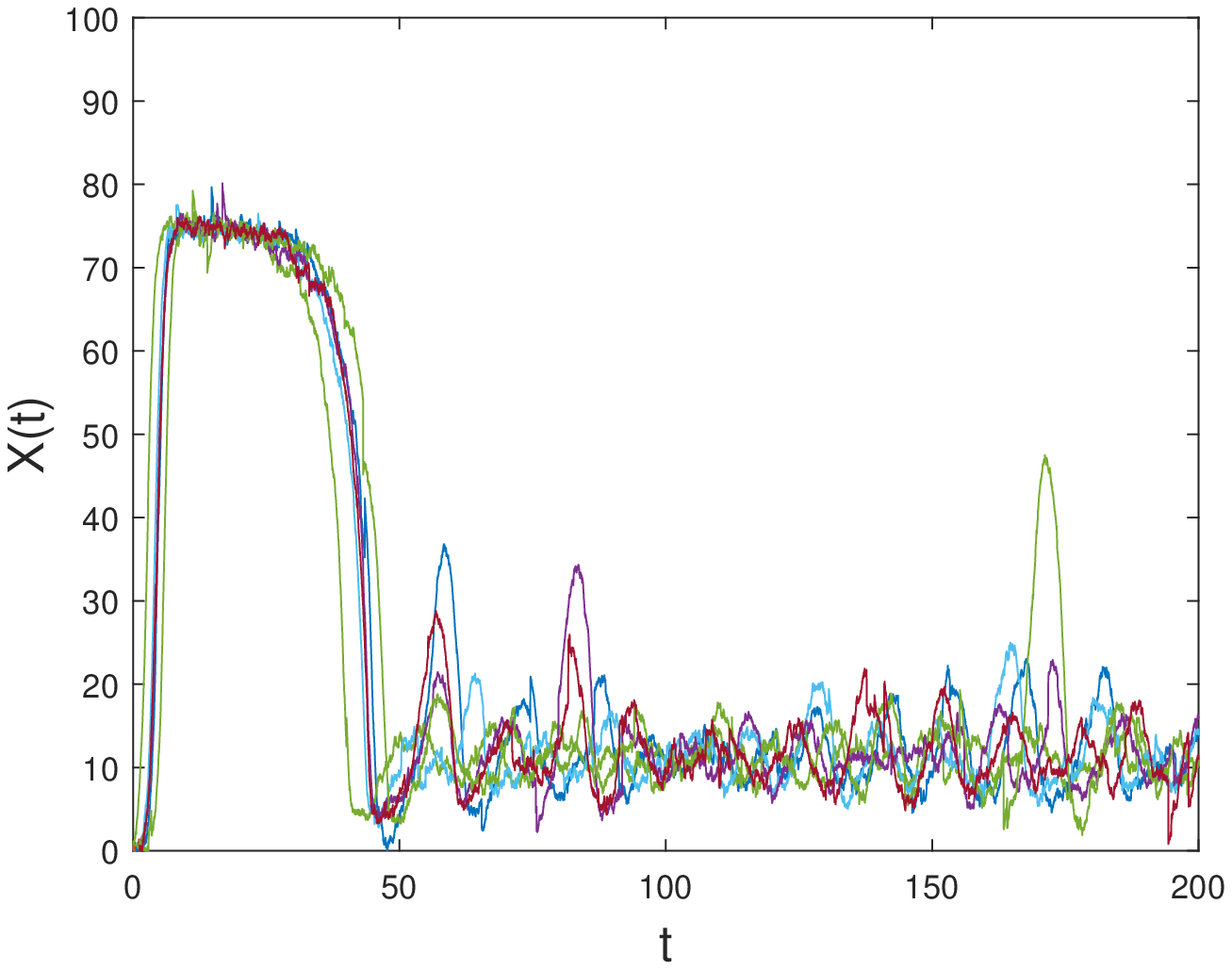}
\end{minipage}
\caption{ The evolution of the prey for a given stability index
$\alpha_1=1.7$ with increasing noise intensity: (a) $\sigma_1=0.04$; (b) $\sigma_1=0.3$; (c) $\sigma_1=0.9$.}\label{Fig10}
\end{center}
\end{figure}

The numerical evolution trajectories of the prey population $X(t)$ in the presence of $\alpha$-stable L\'evy noise $L_{t}^{\alpha_1}$ are shown in Fig \ref{Fig10}. To better understand the effects of the noise intensity $\sigma_1$, we perform simulations by keeping the stability index at constant $\alpha_1=1.7$.  Fig \ref{Fig10}(a) with $\sigma_1=0.04$ displays that the trajectories of $X(t)$ for one set of initial conditions grow to $75$ at the beginning, and then stay at a high level that this circumstance corresponds to the high prey abundance. But some time later, those trajectories present the tendency of decrease. With the increase of time, they vibrate at a gradually declining frequency to reach the low prey abundance around $11.06$ with small-amplitude fluctuations. When $\sigma_1=0.3$, there are some slight bumpiness in the trajectories where
the abundance of prey species is high. But the choppiness is somewhat more intense at the lower levels of prey, which is confirmed numerically in Fig \ref{Fig10}(b). Considering the case for $\sigma_1=0.9$ as in Fig \ref{Fig10}(c), the motion of prey is more vigorous and extensive. Light turbulence happens on the trajectories when the prey is at high abundance. The turbulence is more pronounced when the prey is at low abundance.
More interestingly still, Fig \ref{Fig10} of varying intensity suggests that the prey with low abundance is more vulnerable to environmental changes than that with high abundance.

\begin{figure}[H]
\begin{center}
  \begin{minipage}{2.1in}
\leftline{(a)}
\includegraphics[width=2.1in]{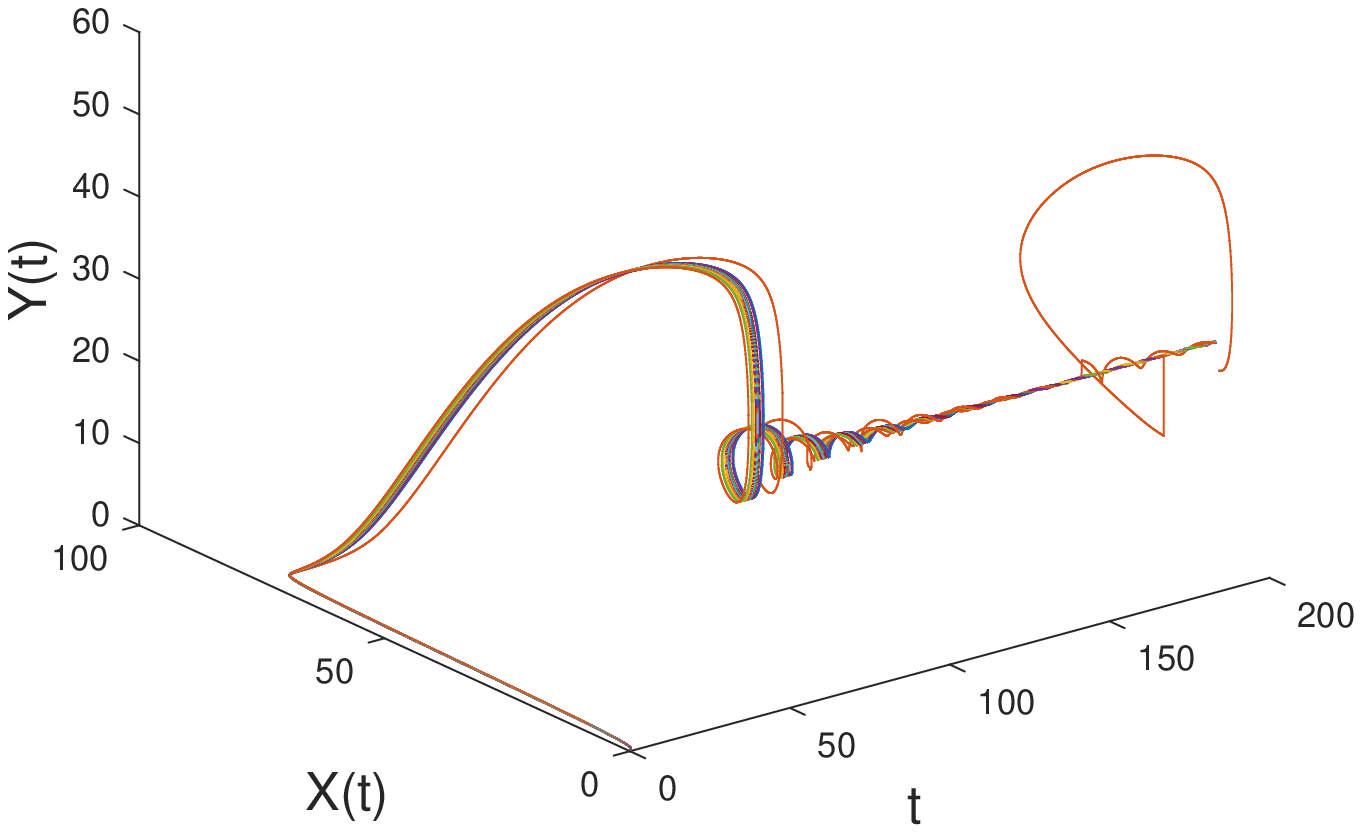}
\end{minipage}
\hfill
\begin{minipage}{2.1in}
\leftline{(b)}
\includegraphics[width=2.1in]{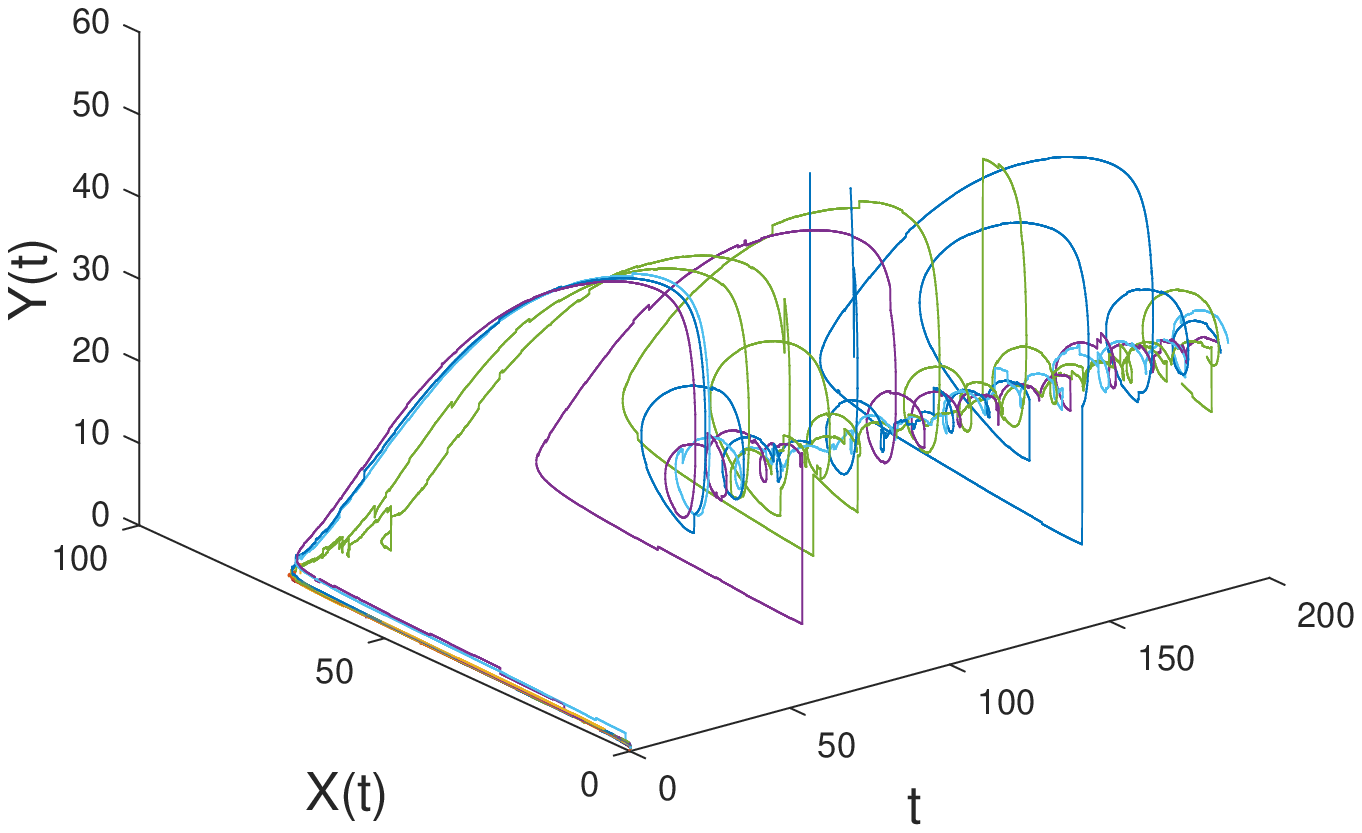}
\end{minipage}
\hfill
  \begin{minipage}{2.1in}
\leftline{(c)}
\includegraphics[width=2.1in]{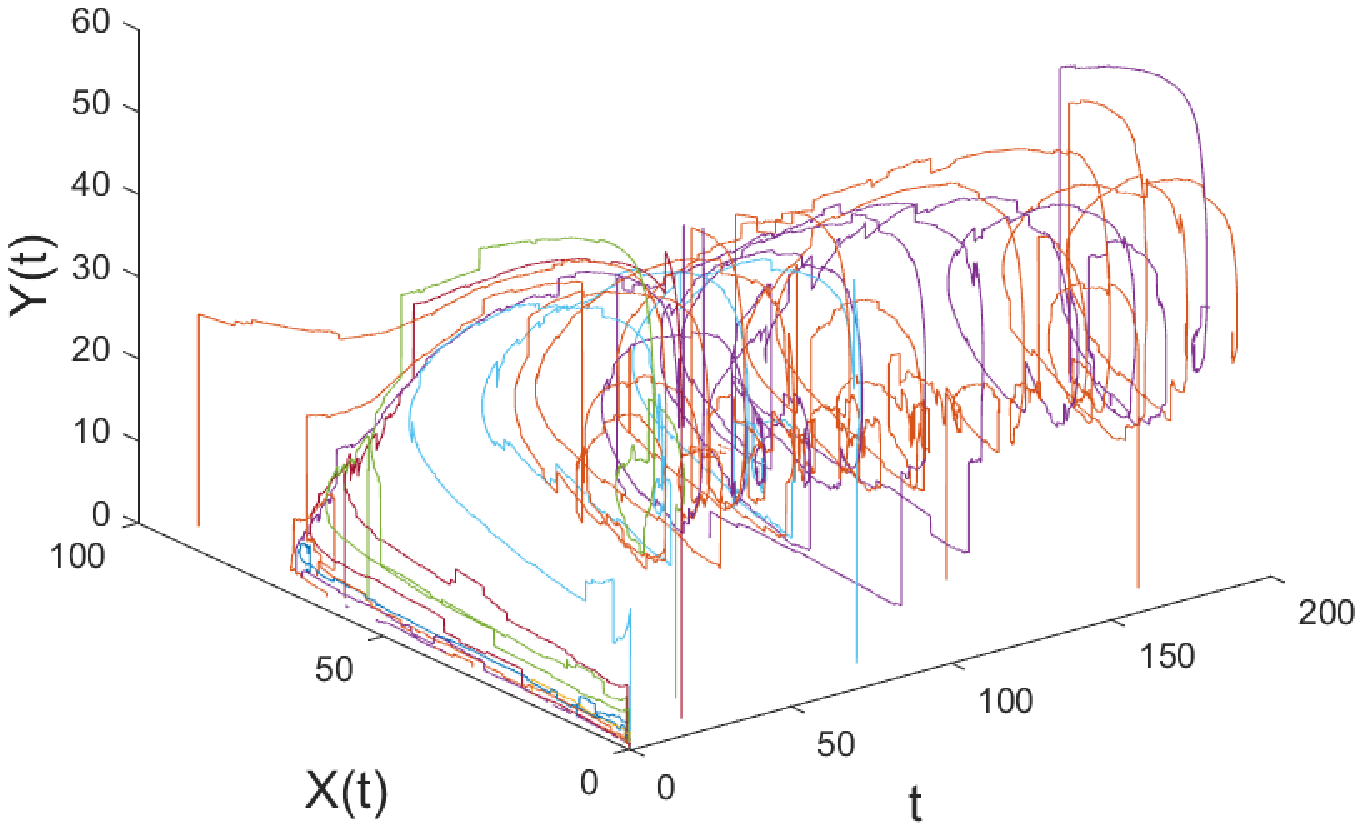}
\end{minipage}
\caption{ The trajectories of two interacting species over time  
for $\sigma_1=0$ (the prey without the perturbation of noise) and $\alpha_2=1$: (a) $\sigma_2=0.001$; (b) $\sigma_2=0.09$; (c) $\sigma_2=0.7$.}\label{Fig11}
\end{center}
\end{figure}

We would like to understand in detail the rich and subtle interplay of the dynamics and the random perturbation  $L_{t}^{\alpha_2}$.
Now we carry out several numerical simulations of stochastic system \eqref{SRM} using the parameter values $\sigma_1=0$ and $\alpha_2=1$. By considering the intensity $\sigma_2=0.001$, Fig \ref{Fig11}(a) gives an indication that a few paths are abnormal. The evolution paths with a set of initial conditions dwell in the vicinity of the $X(t)$-axis at the time of starting,  approach equilibrium point $Z_2$, but then continue to move
towards equilibrium point $Z_3$.  Because of a slight change in the intensity, the noise deteriorates the paths and
leads to significantly different future behavior with respect to $\sigma_2=0.09$, as illustrated in Fig \ref{Fig11}(b). For the parameter value $\sigma_2=0.7$, the system \eqref{SRM} dramatically evolves in a chaotic manner as depicted in Fig \ref{Fig11}(c).

\begin{figure}[H]
\begin{center}
  \begin{minipage}{2.1in}
\leftline{(a)}
\includegraphics[width=2.1in]{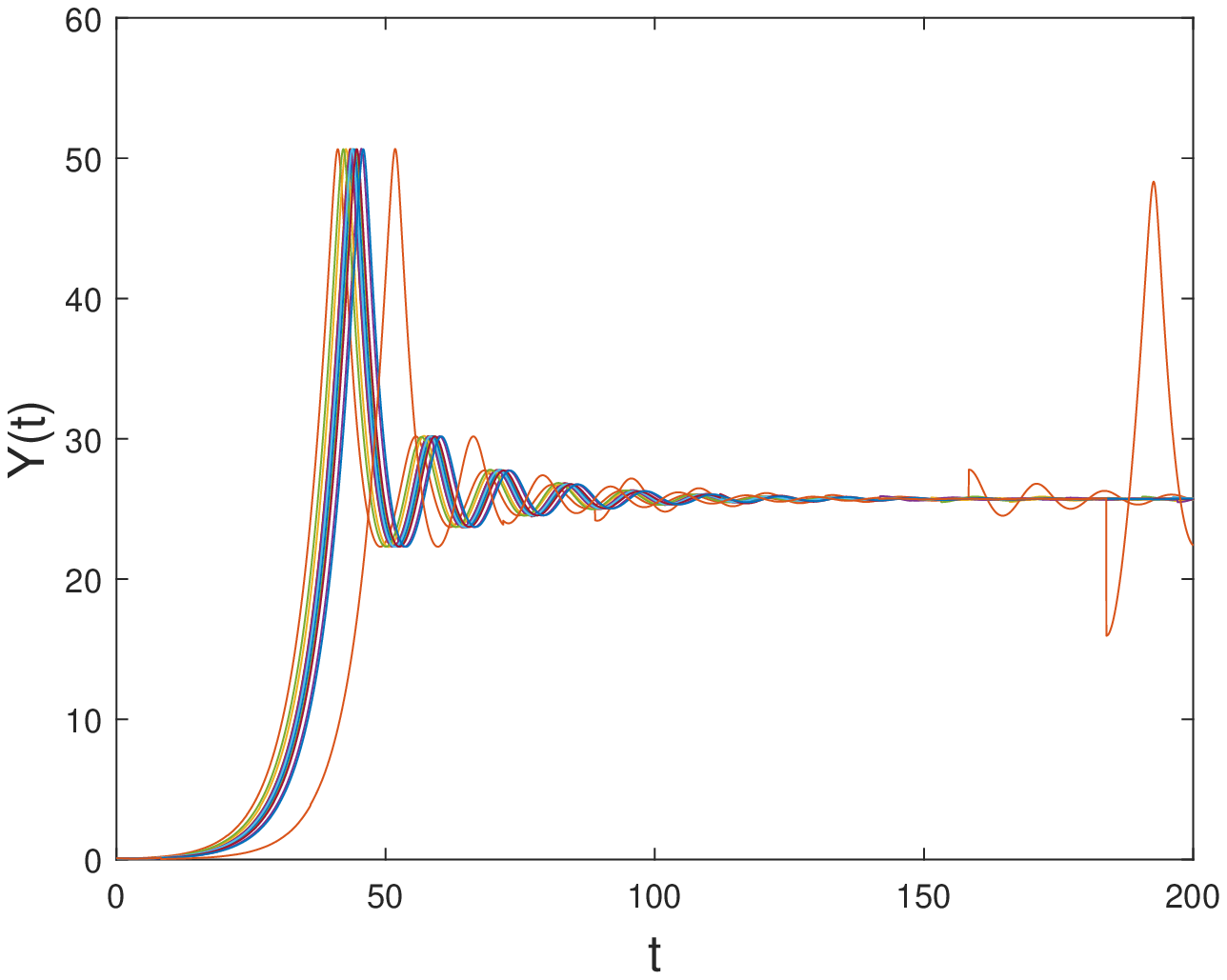}
\end{minipage}
\hfill
\begin{minipage}{2.1in}
\leftline{(b)}
\includegraphics[width=2.1in]{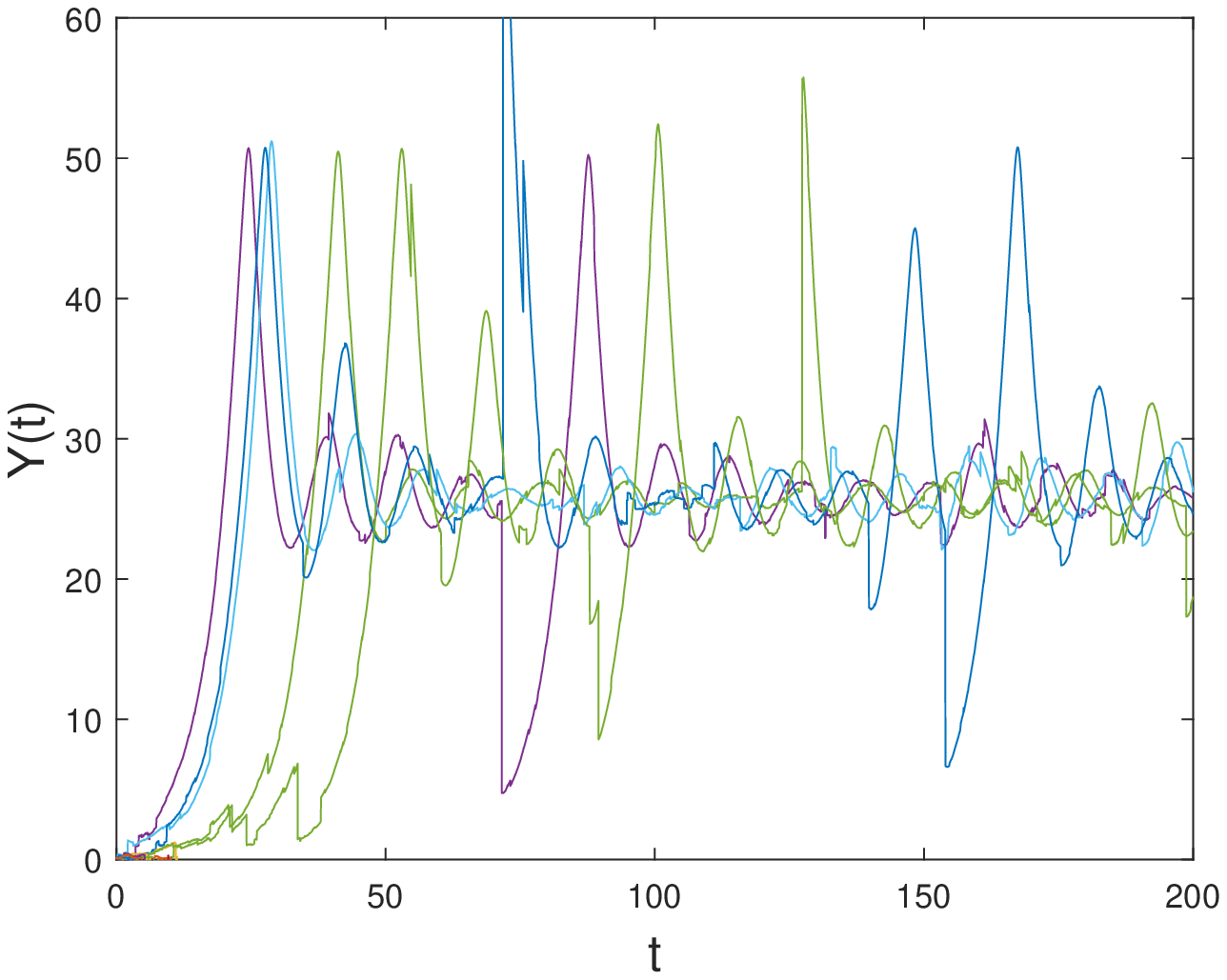}
\end{minipage}
\hfill
  \begin{minipage}{2.1in}
\leftline{(c)}
\includegraphics[width=2.1in]{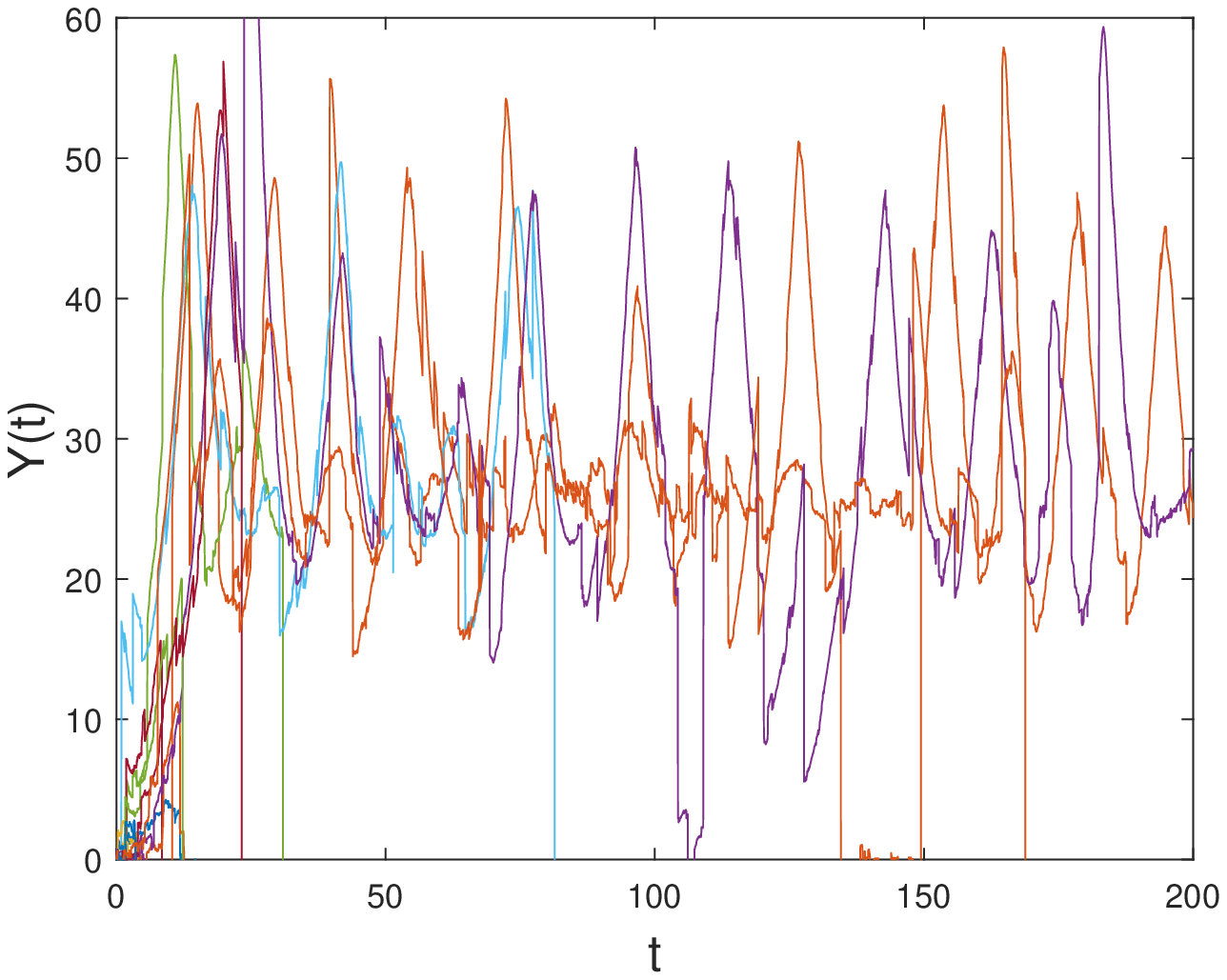}
\end{minipage}
\caption{The evolution of the predator with a fixed stability index
$\alpha_2=1$ as the noise intensity increases: (a) $\sigma_2=0.001$; (b) $\sigma_2=0.09$; (c) $\sigma_2=0.7$.}\label{Fig13}
\end{center}
\end{figure}

By utilizing Monte Carlo method, we perform some dynamical analysis on the predator species as the intensity $\sigma_2$ of $\alpha$-stable L\'evy noise $L_{t}^{\alpha_2}$ varies. Initially, population trajectories for the predator displayed in Fig \ref{Fig13}(a) with $\sigma_2=0.001$ increase fast enough to arrive at 50, but later they decrease rapidly and swing to 25.7. Some of these trajectories are a little naughty, but they are not out of control. As seen in Fig \ref{Fig13}(b), for a higher value of the parameter $\sigma_2=0.09$, the population dynamics of the predator change qualitatively, which are reflected by the apparent randomness of the paths. In Fig \ref{Fig13}(c) we plot the trajectories of the predator species
for $\sigma_2=0.7$. They fluctuate rapidly lacking an ordered organization. This chaotic behavior indicates that the large intensities of L\'evy noises are responsible for large variations in the dynamics.

\subsection{Influence of stability  indexes}
\begin{figure}[H]
\begin{center}
  \begin{minipage}{2.1in}
\leftline{(a)}
\includegraphics[width=2.1in]{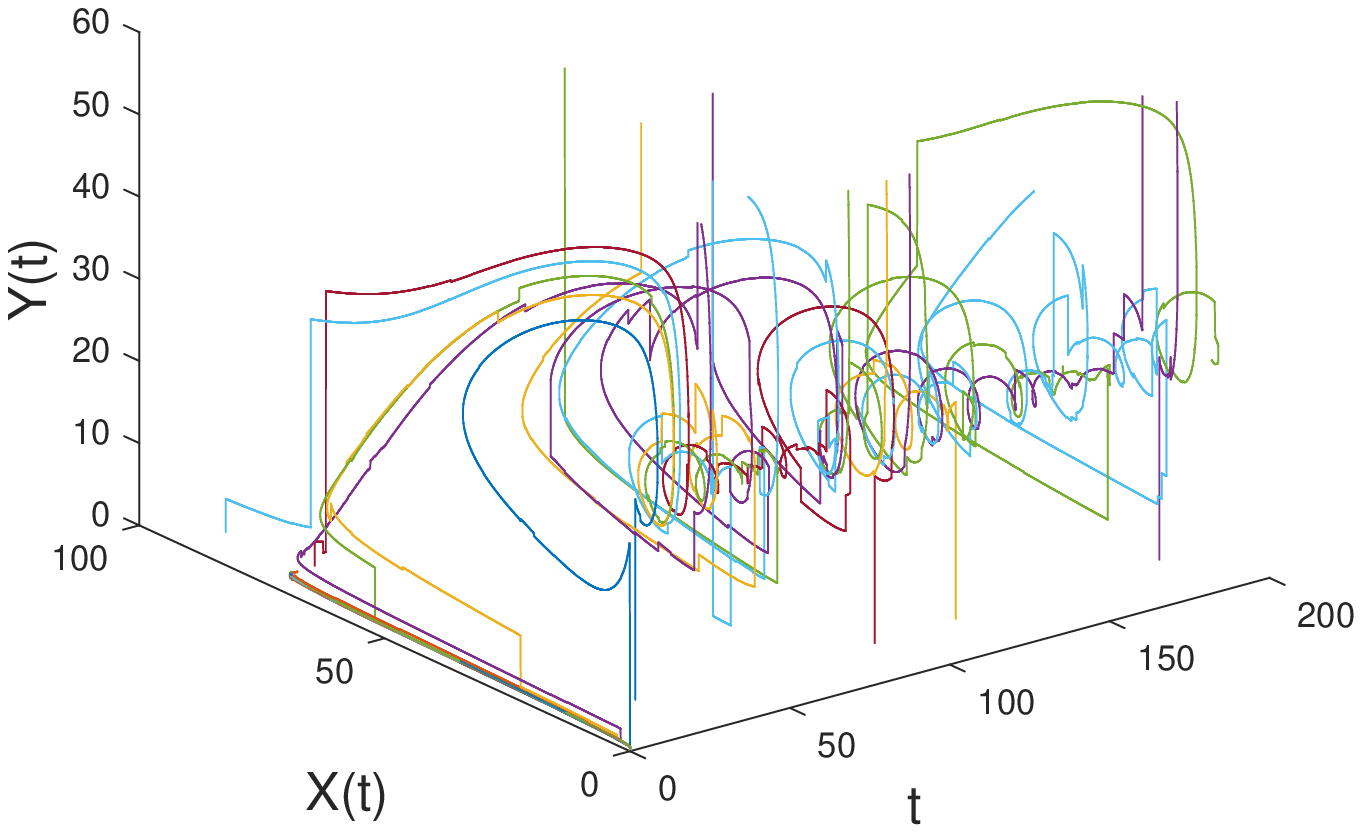}
\end{minipage}
\hfill
\begin{minipage}{2.1in}
\leftline{(b)}
\includegraphics[width=2.1in]{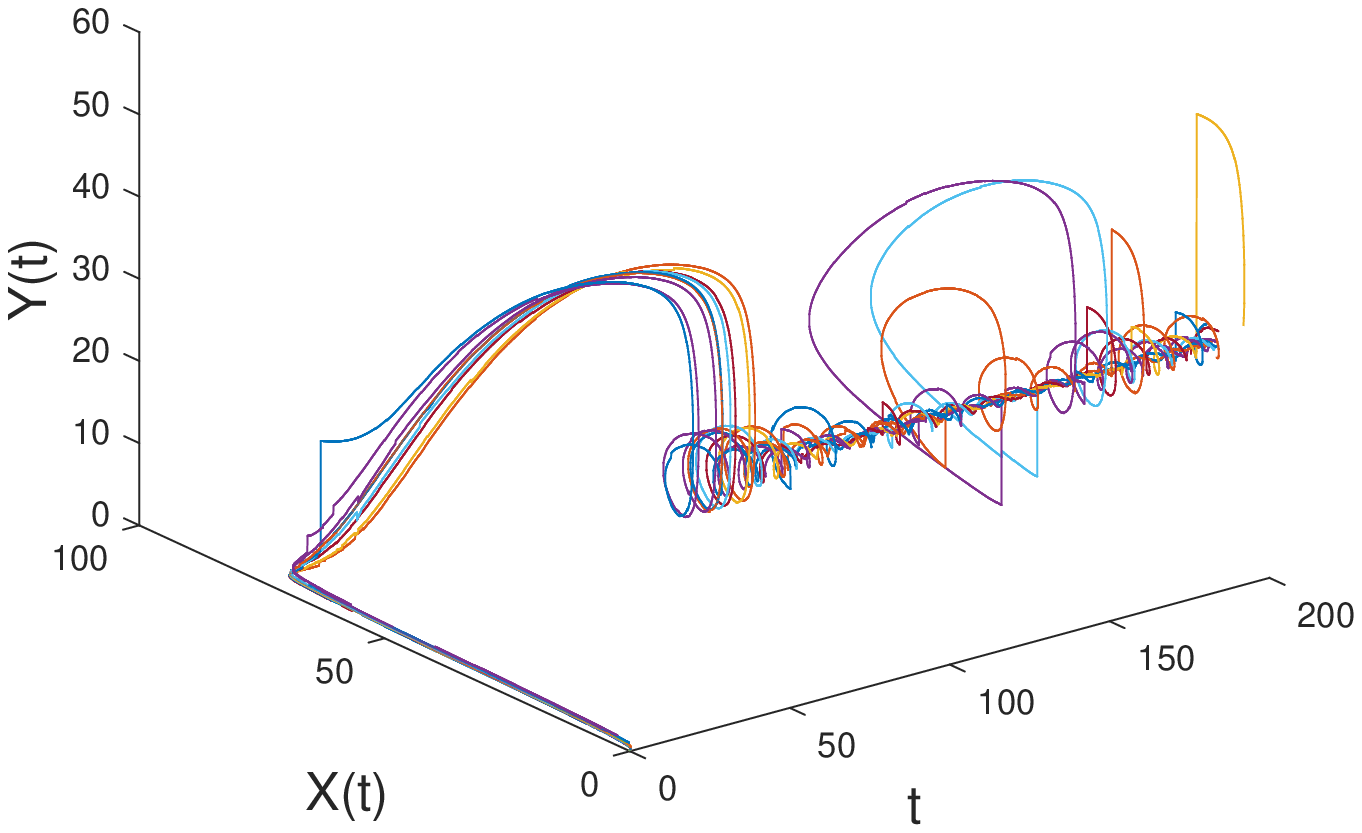}
\end{minipage}
\hfill
  \begin{minipage}{2.1in}
\leftline{(c)}
\includegraphics[width=2.1in]{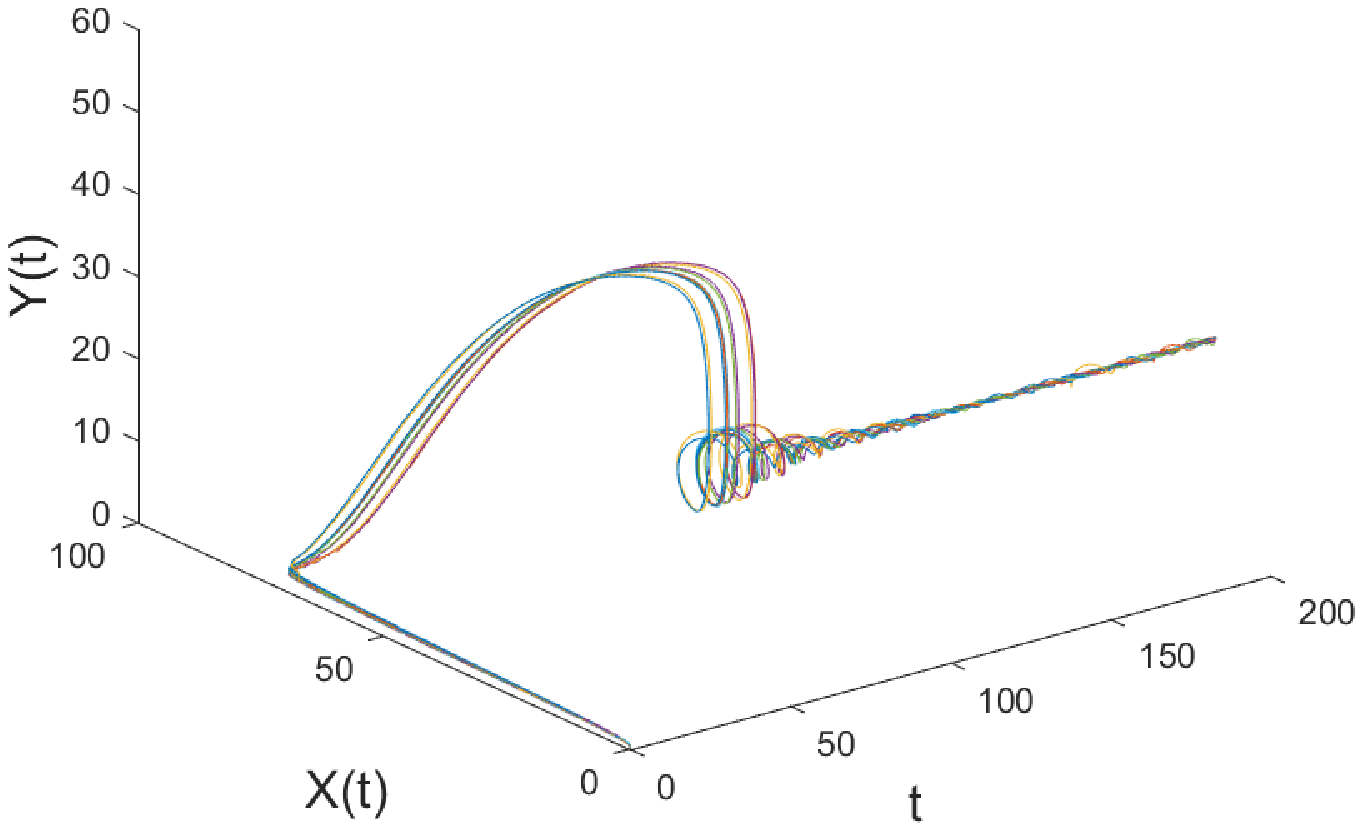}
\end{minipage}
\caption{The trajectories of prey-predator populations over time
with $\sigma_1=0$ and $\sigma_2=0.05$ for the following set of parameter values: (a) $\alpha_2=0.6$; (b) $\alpha_2=1.2$; (c) $\alpha_2=1.8$.}\label{Fig6}
\end{center}
\end{figure}

To show the dynamics of system \eqref{SRM} in which the prey population is noise-free but the predator $Y(t)$ is affected by the noise $L_{t}^{\alpha_2}$, we numerically simulate the curves of two competing species.
In the parameter regime associated with the L\'evy noise intensities $\sigma_1=0$ and $\sigma_2=0.05$, we plot solution curves modeling interacting species. The numerical solutions behave in a complex manner due to the presence of noise. Noise is not applied to the prey $X(t)$ but stochastic disturbance of the predator population can significantly influence the dynamics of the whole prey-predator model. Noise-induced chaotic dynamics are displayed outwardly in Fig \ref{Fig6}(a) for the stability index $\alpha_2=0.6$.
 Small initial differences yield widely diverging outcomes in system \eqref{SRM}. As we can see from Fig \ref{Fig6}(b), the dynamical evolution of two species in competition exhibits somehow complex spatio-temporal oscillations with $\alpha_2=1.2$. The pattern of the populations over time is full of twists and turns. While the opposite behavior occurs for $\alpha_2=1.8$, different trajectories remain close even if they are slightly disturbed. Small initial differences result in small differences of trajectories during a finite time interval demonstrated in Fig \ref{Fig6}(c).

\begin{figure}[H]
\begin{center}
  \begin{minipage}{2.1in}
\leftline{(a)}
\includegraphics[width=2.1in]{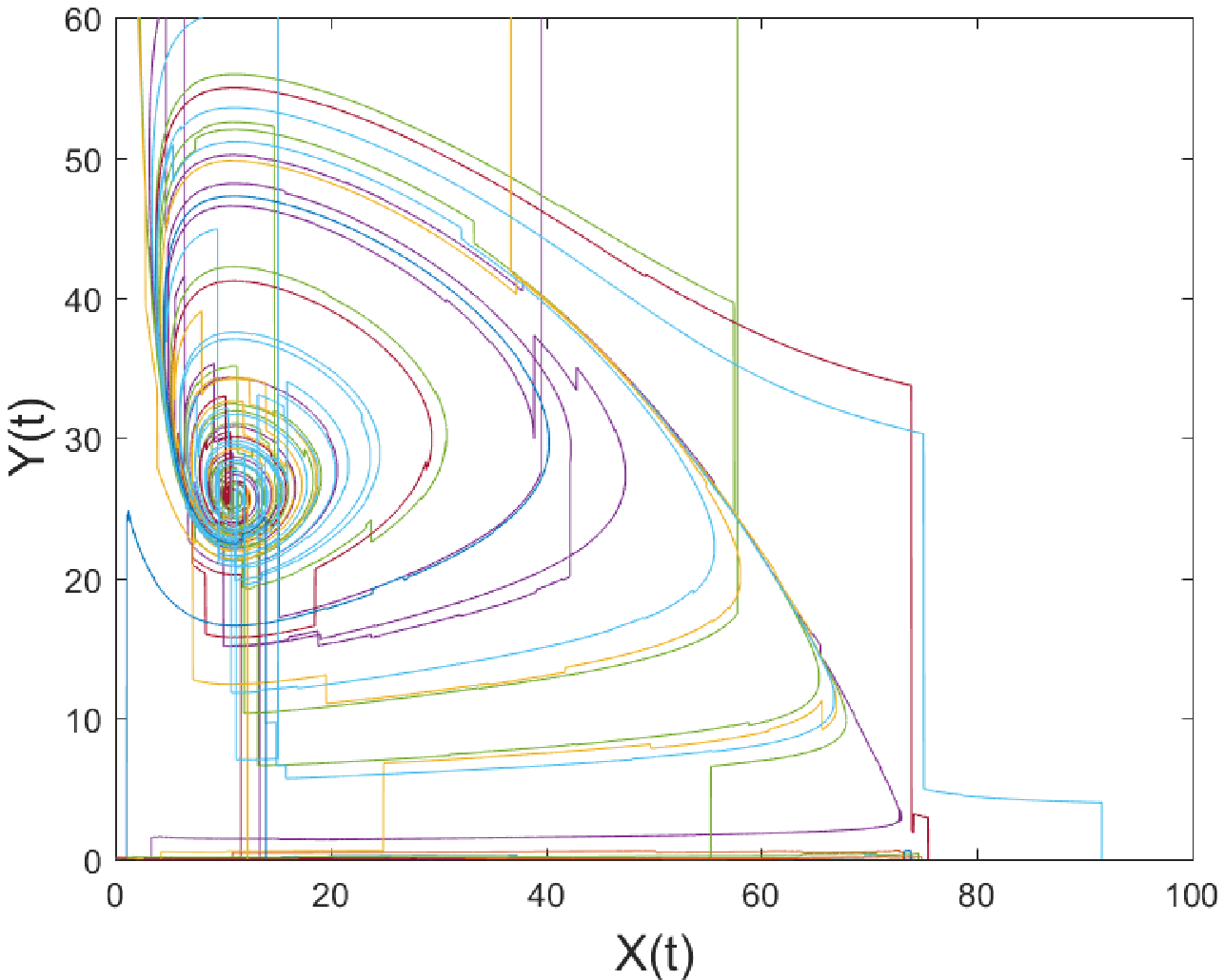}
\end{minipage}
\hfill
\begin{minipage}{2.1in}
\leftline{(b)}
\includegraphics[width=2.1in]{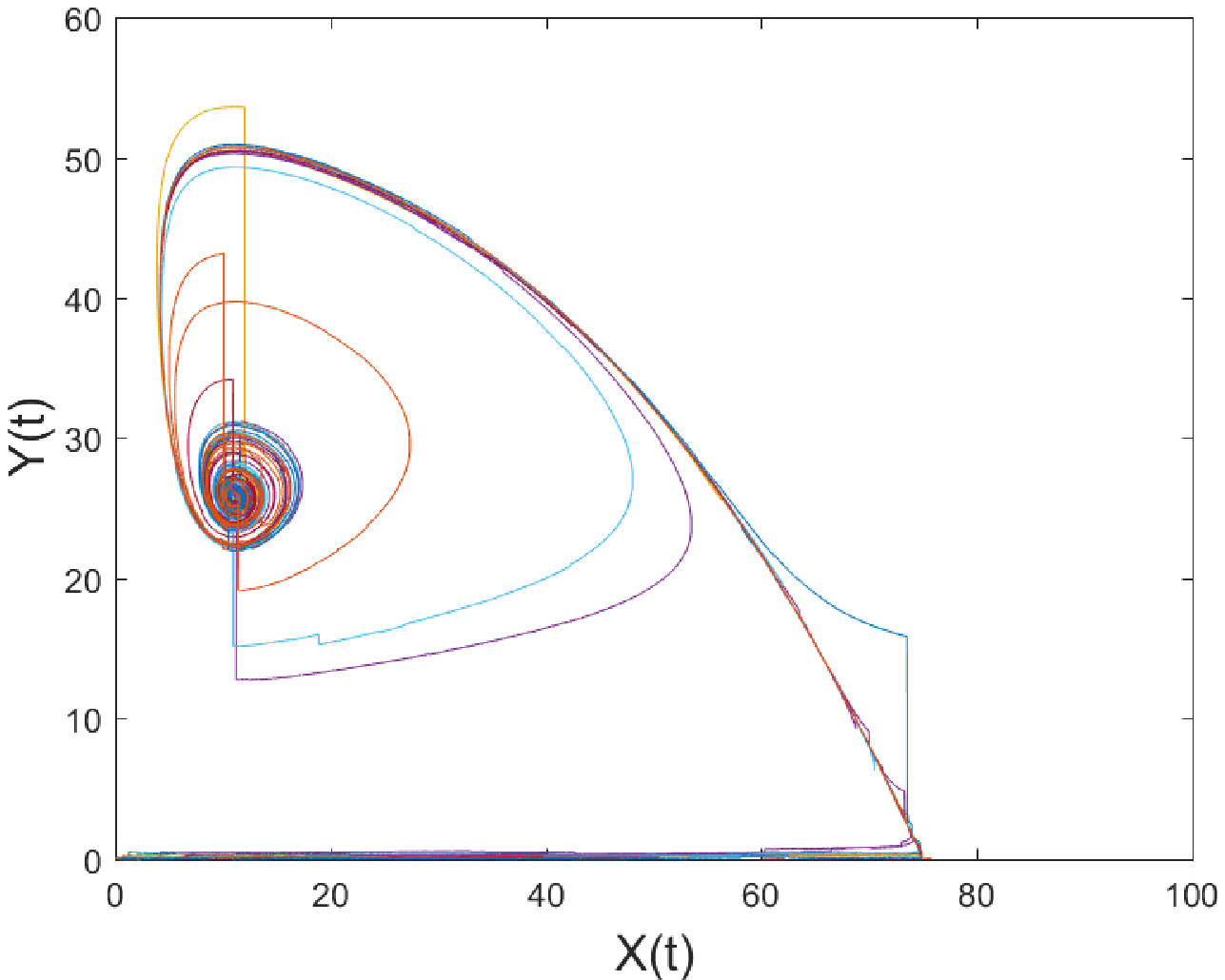}
\end{minipage}
\hfill
  \begin{minipage}{2.1in}
\leftline{(c)}
\includegraphics[width=2.1in]{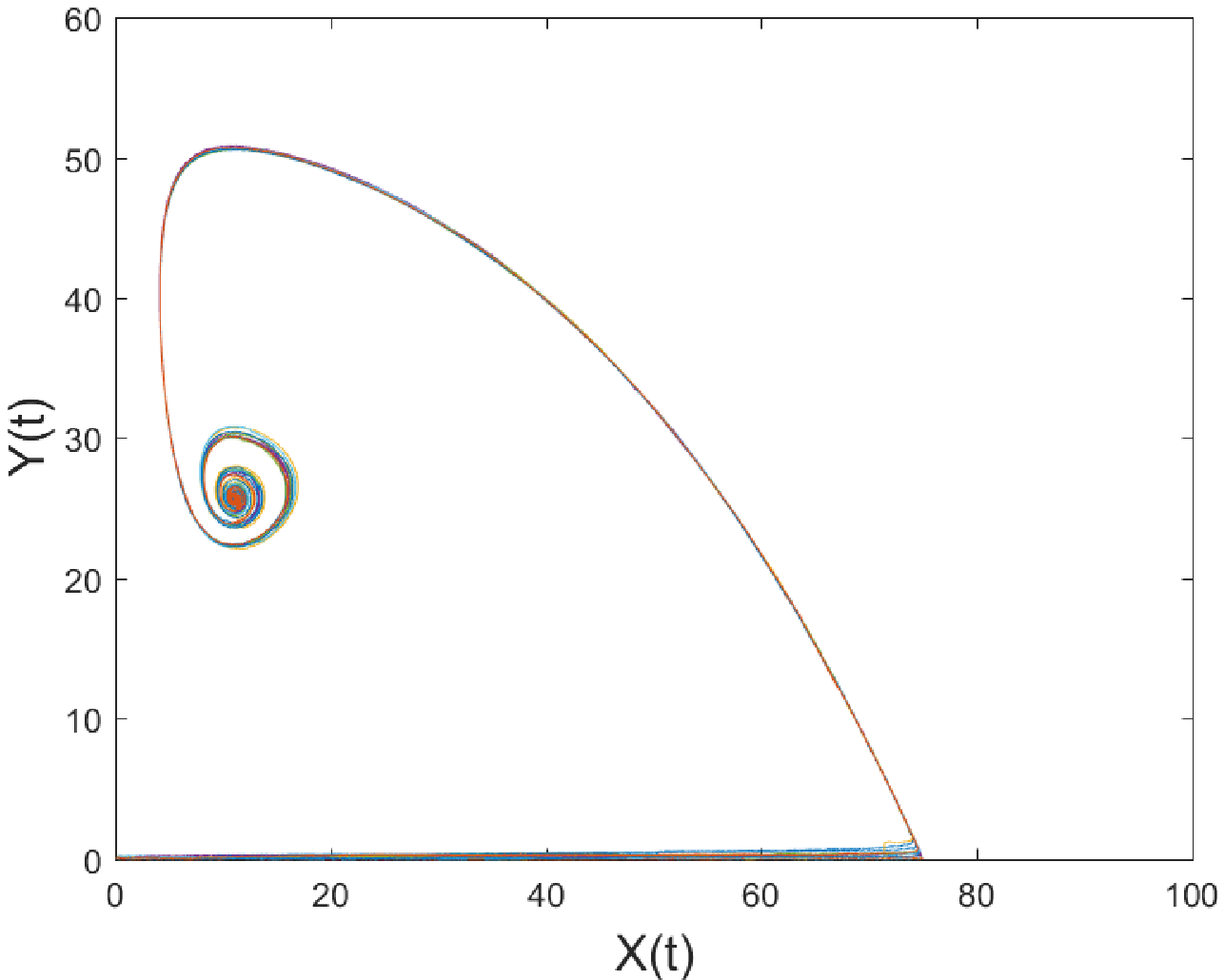}
\end{minipage}
\caption{ The interacting populations $X(t)$ and $Y(t)$
 when $\sigma_1=0$ (the prey in the absence of noise) and $\sigma_2=0.05$: (a) $\alpha_2=0.6$; (b) $\alpha_2=1.2$; (c) $\alpha_2=1.8$.}\label{Fig7}
\end{center}
\end{figure}

Under the same setting as that of Fig \ref{Fig6}, we sketch the prey-predator interactions with parameter values $\sigma_1=0$ and $\sigma_2=0.05$ in Fig \ref{Fig7}. To take into account these interactions,  we choose three different values of $\alpha_2$, i.e.,  $\alpha_2=0.6$, $\alpha_2=1.2$, and $\alpha_2=1.8$. We interpret the results in terms of species and behaviors. In the development of nonlinear stochastic dynamics with the stability index $\alpha_2=0.6$, the creation of chaos by noise is clearly presented in Fig \ref{Fig7}(a), which  portrays that the parameter perturbation can trigger the unordered paths operating in a chaotic regime.
 As depicted in Fig \ref{Fig7}(b), the spiral paths exhibit oscillating patterns with regard to the bigger value $\alpha_2=1.2$. The amazing thing is that oscillating states around $Z_3=(11.06,25.7)$ produce an unusual ear shape. When $\alpha_2=1.8$, the system paths are different from that in the last case. System curves with starting conditions near the origin equilibrium $Z_1$ are horizontal until they arrive at $Z_2=(75,0)$. After that, those curves tend to evolve toward $Z_3=(11.06,25.7)$ along the spiral, as clearly detailed in Fig \ref{Fig7}(c). For the increasing $\alpha_2$, our computations reveal that changes happen abruptly.

\begin{figure}[H]
\begin{center}
  \begin{minipage}{2.1in}
\leftline{(a)}
\includegraphics[width=2.1in]{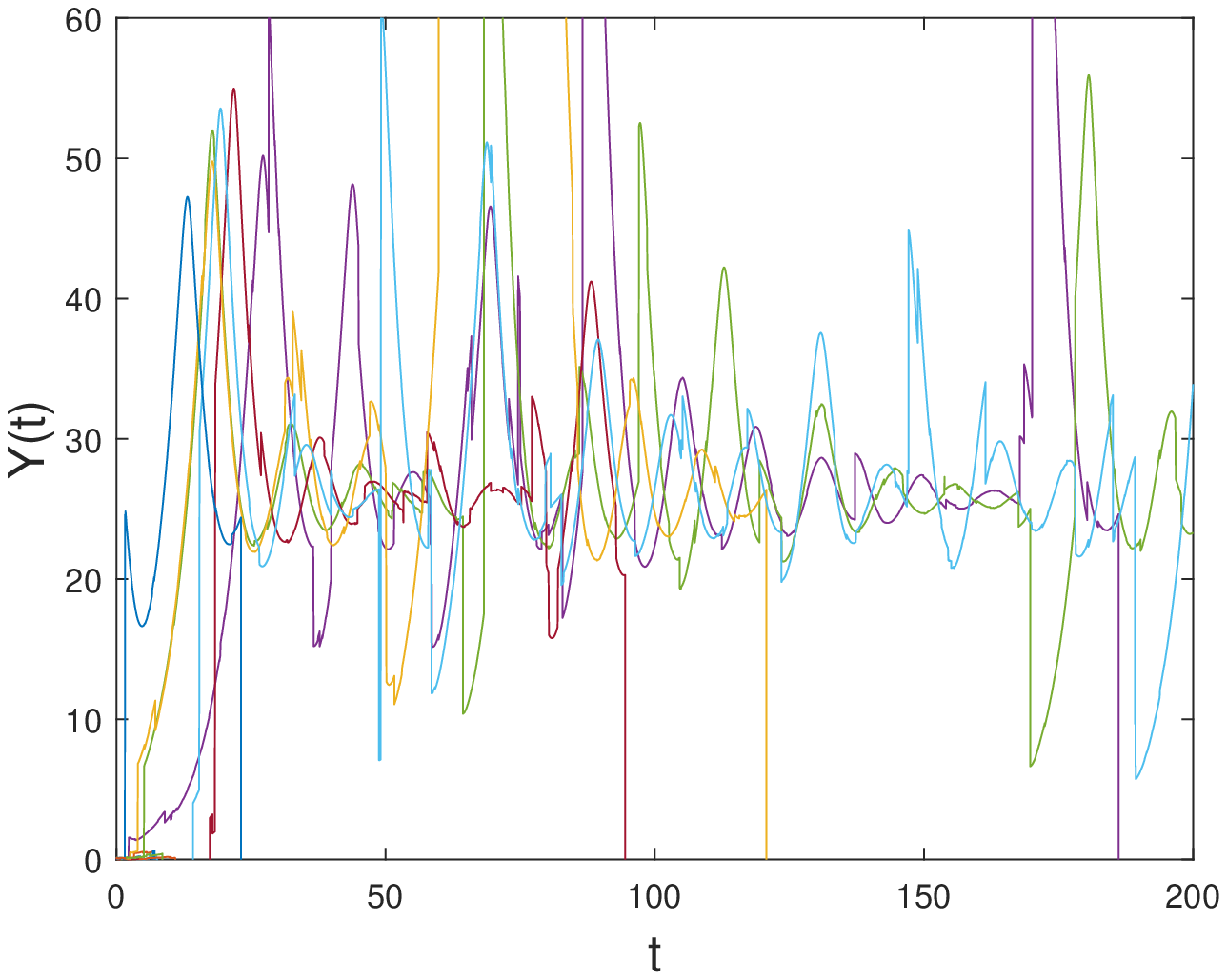}
\end{minipage}
\hfill
\begin{minipage}{2.1in}
\leftline{(b)}
\includegraphics[width=2.1in]{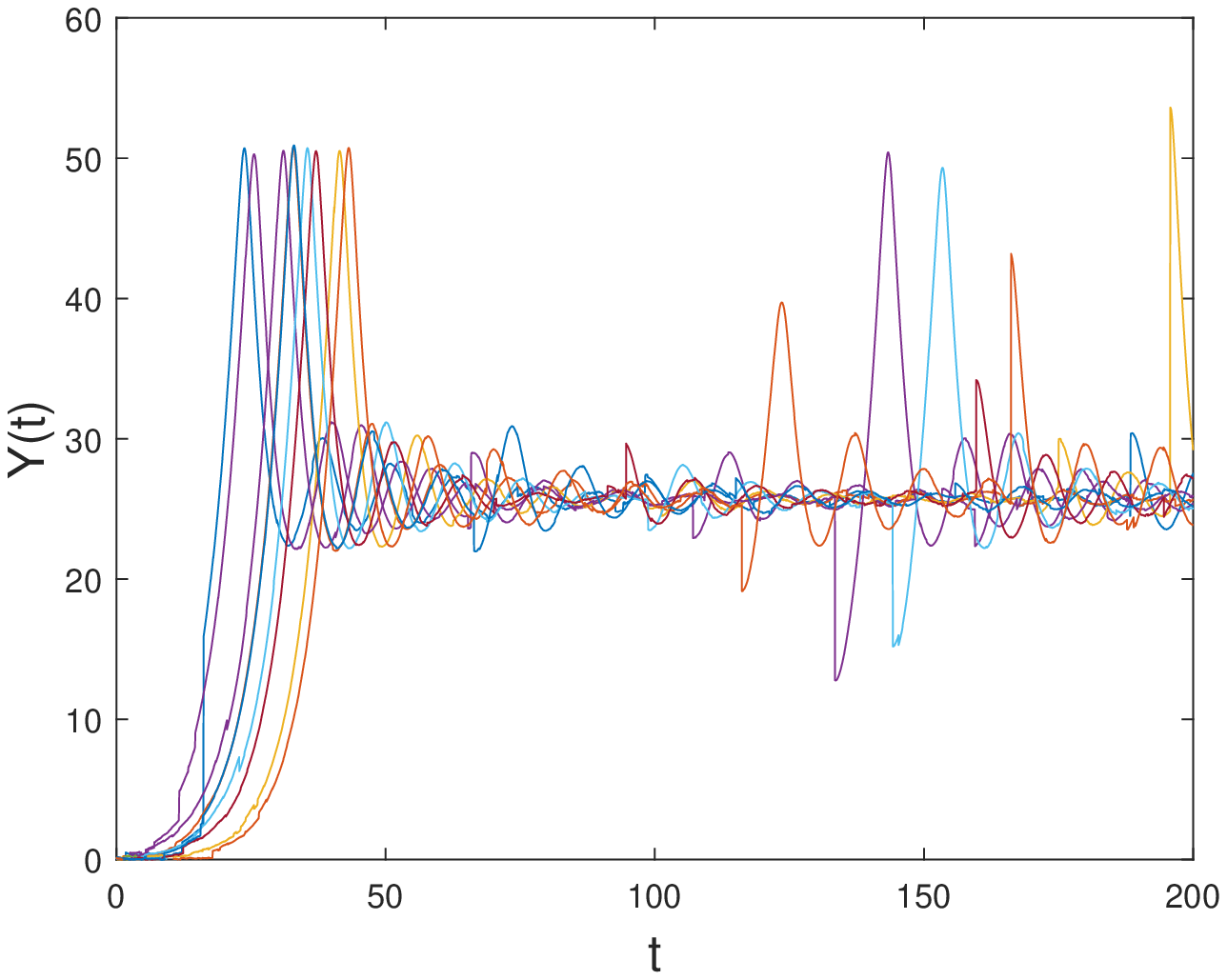}
\end{minipage}
\hfill
  \begin{minipage}{2.1in}
\leftline{(c)}
\includegraphics[width=2.1in]{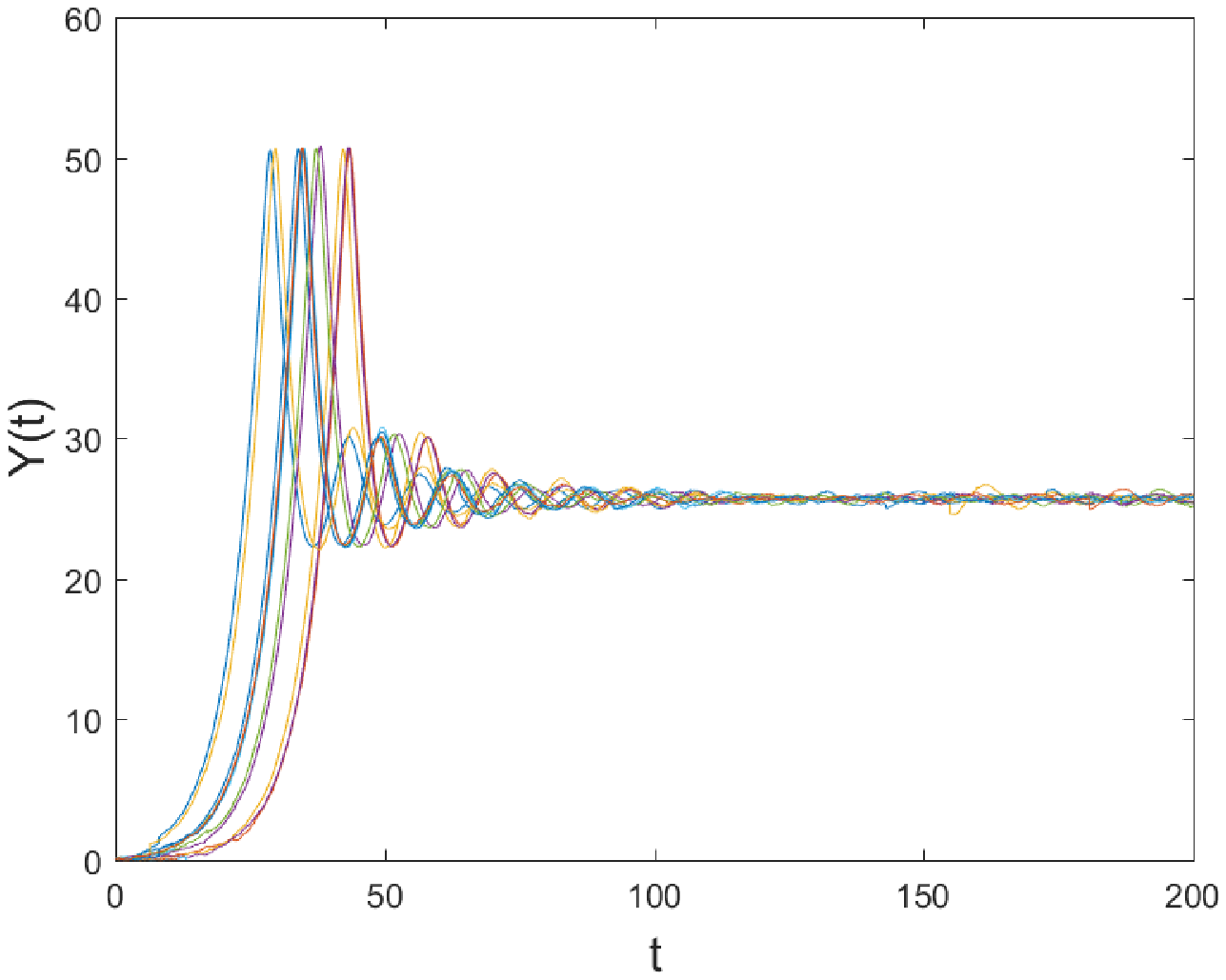}
\end{minipage}
\caption{The evolution of the predator at $\sigma_2=0.05$ with a gradual increase in the stability index: (a) $\alpha_2=0.6$; (b) $\alpha_2=1.2$; (c) $\alpha_2=1.8$.}\label{Fig8}
\end{center}
\end{figure}

Specifically, we explore the complicated dynamics of the predator population $Y(t)$ contaminated by noise by examining a wide variety of paths. We still choose the noise intensity $\sigma_2=0.05$. We observe from Fig \ref{Fig8}(a) that for $\alpha_2=0.6$ the paths have jumps with higher frequencies, which describe chaotic behaviours. We also find that the noise $L_{t}^{\alpha_2}$ blurs deterministic solutions of $Y(t)$ when $\alpha_2=1.2$ as in Fig \ref{Fig8}(b).  While for $\alpha_2=1.8$, all trajectories of $Y(t)$ with initial conditions at $t=0$ grow logistically approaching $50$ in the early stage, and then they are reduced to 25.7 accompanying fluctuations with lower probabilities; see Fig \ref{Fig8}(c).

\section{Conclusions and future challenges}\label{Cf}
We successfully explored a Rosenzweig-MacArthur prey-predator model. It was discussed that the proposed model \eqref{RM} has at most three equilibrium points, i.e., the extinction of population $Z_1$, the predator-free point $Z_2$, and the coexistence point $Z_3$. The equilibrium points $Z_2$ and $Z_3$ are conditionally asymptotically stable. Our analysis also showed that the model \eqref{RM} may undergo a transcritical bifurcation for suitable parameter values. Further, we still obtained the non-local Fokker-Planck equation analytically, but we used a numerical integration method to get a complete view about the stationary density. More importantly, we investigated the chaotic dynamics of stochastic model \eqref{SRM} and provided insights into the effect of varying the values of the parameters. We carried out several numerical simulations to corroborate our findings.

Our analytical findings and numerical simulations can be extended to other disciplines to address non-Gaussian stochastic systems. In fact, the methods presented here are rather general and can also be used to work on population models for interacting species with other and more general nonlinearities. The noises may drastically modify the deterministic dynamics, but we only focus on finding that $\alpha$-stable L\'evy noises induce chaotic behavior. A further analysis of this would indeed be worthwhile. Chaotic systems share many properties with noisy systems, which could be of independent interest. In reality, the prey-predator dynamical systems experience influence from their stochastic environments. Motivated by these evidentiary statistics, we have to make some changes to protect the environment and prevent environmental deterioration now and in the future.

\section*{Appendix}

Applying It\^o's formula to stochastic dynamical system \eqref{SRM}, we establish
{\small\begin{align}\nonumber
df(X_{t},Y_{t})
&=(rX_{t}-0.02X_{t}^{2}-Y_{t}f(X_{t}))\frac{\partial}{\partial X_{t}}f(X_{t},Y_{t})dt+(0.4f(X_{t})-\mu)Y_{t}\frac{\partial}{\partial Y_{t}}f(X_{t},Y_{t})dt\\ \nonumber
&\,\,\,\,\,\,+\int_{\mathbb{R}\setminus \{ 0 \}}(f(X_{t}+\sigma_{1}u,Y_{t})-f(X_{t},Y_{t})-\sigma_{1}u\chi_{\{|u|\leq1\}}\frac{\partial}{\partial X_{t}}f(X_{t},Y_{t}))\nu_{\alpha_1}(du)dt\\ \nonumber
&\,\,\,\,\,\,+\int_{\mathbb{R}\setminus \{ 0 \}}(f(X_{t},Y_{t}+\sigma_{2}v)-f(X_{t},Y_{t})-\sigma_{2}v\chi_{\{|v|\leq1\}}\frac{\partial}{\partial Y_{t}}f(X_{t},Y_{t}))\nu_{\alpha_2}(dv)dt\\ \nonumber
&=[(rX_{t}-0.02X_{t}^{2}-Y_{t}f(X_{t}))\frac{\partial}{\partial X_{t}}f(X_{t},Y_{t})+(0.4f(X_{t})-\mu)Y_{t}\frac{\partial}{\partial Y_{t}}f(X_{t},Y_{t})\\  \nonumber
&\,\,\,\,\,\,+\sigma_1^{\alpha_1}\int_{\mathbb{R}\setminus \{ 0 \}}(f(X_{t}+u,Y_{t})-f(X_{t},Y_{t}))\nu_{\alpha_1}(du)\\ \label{free}
&\,\,\,\,\,\,+\sigma_2^{\alpha_2}\int_{\mathbb{R}\setminus \{ 0 \}}(f(X_{t},Y_{t}+v)-f(X_{t},Y_{t}))\nu_{\alpha_2}(dv)]dt,
\end{align}}
where $\chi_{\{|u|\leq1\}}$ is the indicator function of the set $\{|u|\leq1\}$. Taking expectation on both sides of (\ref{free}), we get
\begin{align}\nonumber
&d\mathbb{E}f(X_{t},Y_{t})\\ \nonumber
&=\mathbb{E}[(rX_{t}-0.02X_{t}^{2}-Y_{t}f(X_{t}))\frac{\partial}{\partial X_{t}}f(X_{t},Y_{t})+(0.4f(X_{t})-\mu)Y_{t}\frac{\partial}{\partial Y_{t}}f(X_{t},Y_{t})\\  \nonumber
&\,\,\,\,\,\,+\sigma_1^{\alpha_1}\int_{\mathbb{R}\setminus \{ 0 \}}(f(X_{t}+u,Y_{t})-f(X_{t},Y_{t}))\nu_{\alpha_1}(du)\\  \label{goon}
&\,\,\,\,\,\,+\sigma_2^{\alpha_2}\int_{\mathbb{R}\setminus \{ 0 \}}(f(X_{t},Y_{t}+v)-f(X_{t},Y_{t}))\nu_{\alpha_2}(dv)]dt.
\end{align}
It is relevant to point out that the generator for system (\ref{SRM}) is
\begin{align*}
  Ap(x,y,t)&=(rx-0.02x^{2}-yf(x))\frac{\partial}{\partial x}p(x,y,t)+(0.4f(x)-\mu)y\frac{\partial}{\partial y}p(x,y,t)\\
           &\,\,\,\,\,\,+\sigma_1^{\alpha_1}\int_{\mathbb{R}\setminus \{ 0 \}}(p(x+u,y,t)-p(x,y,t))\nu_{\alpha_1}(du)\\
           &\,\,\,\,\,\,+\sigma_2^{\alpha_2}\int_{\mathbb{R}\setminus \{ 0 \}}(p(x,y+v,t)-p(x,y,t))\nu_{\alpha_2}(dv).
\end{align*}
We rewrite the equation (\ref{goon}) into
\begin{align}\nonumber
&\frac{d}{dt}\mathbb{E}f(X_{t}^{\varepsilon},Y_{t}^{\varepsilon})\\ \nonumber
&=\mathbb{E}[(rX_{t}-0.02X_{t}^{2}-Y_{t}f(X_{t}))\frac{\partial}{\partial X_{t}}f(X_{t},Y_{t})+(0.4f(X_{t})-\mu)Y_{t}\frac{\partial}{\partial Y_{t}}f(X_{t},Y_{t})\\  \nonumber
&\,\,\,\,\,\,+\sigma_1^{\alpha_1}\int_{\mathbb{R}\setminus \{ 0 \}}(f(X_{t}+u,Y_{t})-f(X_{t},Y_{t}))\nu_{\alpha_1}(du)\\ \nonumber
&\,\,\,\,\,\,+\sigma_2^{\alpha_2}\int_{\mathbb{R}\setminus \{ 0 \}}(f(X_{t},Y_{t}+v)-f(X_{t},Y_{t}))\nu_{\alpha_2}(dv)]\\ \nonumber
&=\int_{\mathbb{R}^{2}}[(rx-0.02x^{2}-yf(x))\frac{\partial}{\partial x}f(x,y)+(0.4f(x)-\mu)y\frac{\partial}{\partial y}f(x,y)\\ \nonumber
&\,\,\,\,\,\,+\sigma_1^{\alpha_1}\int_{\mathbb{R}\setminus \{ 0 \}}(f(x+u,y)-f(x,y))\nu_{\alpha_1}(du)\\ \nonumber
&\,\,\,\,\,\,+\sigma_2^{\alpha_2}\int_{\mathbb{R}\setminus \{ 0 \}}(f(x,y+v)-f(x,y))\nu_{\alpha_2}(dv)]p(x,y,t)dxdy\\ \nonumber
&=\int_{\mathbb{R}^{2}}f(x,y)\Big(-\frac{\partial}{\partial x}[(rx-0.02x^{2}-yf(x))p(x,y,t)]-\frac{\partial}{\partial y}[(0.4f(x)-\mu)yp(x,y,t)]\\ \nonumber
&\,\,\,\,\,\,-\sigma_1^{\alpha_1}\int_{\mathbb{R}\setminus \{ 0 \}}(p(x,y,t)-p(x-u,y,t))\nu_{\alpha_1}(du)\\ \nonumber
&\,\,\,\,\,\,-\sigma_2^{\alpha_2}\int_{\mathbb{R}\setminus \{ 0 \}}(p(x,y,t)-p(x,y-v,t))\nu_{\alpha_2}(dv)\Big)dxdy.
\end{align}
Observe that the adjoint operator of the generator $A$ is
\begin{align*}
  A^{*}p(x,y,t)&=-\frac{\partial}{\partial x}[(rx-0.02x^{2}-yf(x))p(x,y,t)]-\frac{\partial}{\partial y}[(0.4f(x)-\mu)yp(x,y,t)]\\ \nonumber
&\,\,\,\,\,\,-\sigma_1^{\alpha_1}\int_{\mathbb{R}\setminus \{ 0 \}}(p(x,y,t)-p(x-u,y,t))\nu_{\alpha_1}(du)\\ \nonumber
&\,\,\,\,\,\,-\sigma_2^{\alpha_2}\int_{\mathbb{R}\setminus \{ 0 \}}(p(x,y,t)-p(x,y-v,t))\nu_{\alpha_2}(dv)\\ \nonumber
&=\big(0.04x+yf'(x)+\mu-0.4f(x)-r\big)p(x,y,t)  \\ \nonumber
	&~~~~+\big(0.02x^{2}+yf(x)-rx\big)\frac{\partial}{\partial x}p(x,y,t)+y(\mu-0.4f(x))\frac{\partial}{\partial y}p(x,y,t)  \\ \nonumber
	&~~~~+\sigma_1^{\alpha_1}\int_{\mathbb{R}\backslash\{0\}}(p(x+u,y,t)-p(x,y,t))\nu_{\alpha_1}(du)\\ \nonumber
&~~~~+\sigma_2^{\alpha_2}\int_{\mathbb{R}\backslash\{0\}}(p(x,y+v,t)-p(x,y,t))\nu_{\alpha_2}(dv).
\end{align*}
Therefore, the Fokker-Planck equation for system (\ref{SRM}) is the equation \eqref{FP}.

\bigskip
\noindent\textbf{DATA AVAILABILITY}

Numerical algorithms source code that support the findings of this study are openly
available in GitHub, Ref. \cite{SY}.

\bigskip
\noindent\textbf{ACKNOWLEDGMENTS}

The authors are happy to thank Jinqiao Duan, Haitao Xu and Zhigang Zeng for fruitful discussions on stochastic dynamical systems. The authors acknowledge support from the NSFC grant 12001213.  







\section*{References}

\bibliography{mybibfile}

\end{document}